\documentclass[11pt]{amsproc}
\newtheorem{thm}{Th\'eor\`eme}[section]
\newtheorem{prop}{Proposition}[section]
\newtheorem{cor}{Corollaire}[section]
\newtheorem{rem}{Remarque}[section]
\newtheorem{lem}{Lemme}[section]
\newtheorem{ex}{Exemple}[section]

\newcommand{\cqfd}
{%
\fbox{\rule[+.2mm]{0.2mm}{0cm}}

\smallskip
}

\numberwithin{equation}{section}

\def\R{{\mathbb R}}

\def\Z{{\mathbb Z}} 
\def\C{{\mathbb C}} 
\def\h{{\zeta \hskip0.2mm}}
\def\d{{\delta \,}}

\begin{document}

\title[Asymptotiques de nombres de Betti]{Asymptotiques de nombres de
  Betti d'hypersurfaces projectives r\'eelles} 
 
\author{F. Bihan}


\email{
Frederic.Bihan@math.unil.ch
}



\begin{abstract} 
On s'int\'eresse \`a la valeur maximale des nombres de Betti
$b_i({\R}X_m^n)$ pour $i,m,n$ fix\'es o\`u ${\R}X_m^n$ est la partie r\'eelle
d'une hypersurface r\'eelle $X_m^n$ non singuli\`ere de degr\'e $m$
dans l'espace projectif complexe ${\C}P^n$, ainsi qu'\`a la valeur maximale
des nombres de Betti
$b_i({\R}Y_{2k}^n)$ pour $i,k,n$ fix\'es o\`u ${\R}Y_{2k}^n$ est la partie r\'eelle
d'un rev\^etement double r\'eel de ${\C}P^n$ ramifi\'e sur une hypersurface
r\'eelle non singuli\`ere de degr\'e $2k$.
On montre l'existence de limites $\lim_{m \rightarrow +\infty} \frac{Max\,  b_i({\R}X_m^n)}{m^n}={\h}_{i,n}$
et $\lim_{k \rightarrow +\infty} \frac{Max\,  b_i({\R}Y_{2k}^n)}{k^n}={\d}_{i,n}.$
On construit des hypersurfaces par petites perturbations
d'hypersurfaces doubles en utilisant la m\'ethode de Viro. Cette construction permet d'obtenir
des bornes inf\'erieures r\'ecursives pour les ${\h}_{0,n}$ et ${\d}_{0,n}$, ainsi que les in\'egalit\'es
$ {\h}_{0,3} \geq \frac{{\d}_{0,2}}{6}+\frac{1}{12}$ et ${\h}_{1,3} \geq \frac{{\d}_{1,2}}{6}+\frac{1}{6}$
concernant les surfaces alg\'ebriques dans ${\C}P^3$. On montre alors que, pour tout $n \geq 5$, il existe
des hypersurfaces r\'eelles $X_m^n \subset {\C}P^n$ qui ne sont pas des $T$-hypersurfaces, et l'on obtient les in\'egalit\'es
$ \frac{35}{96} \leq {\h}_{0,3} \leq \frac{5}{12}$ et $\frac{35}{48} \leq {\h}_{1,3} \leq \frac{5}{6}$.
\end{abstract}

\maketitle
\specialsection*{Introduction}

Une hypersurface alg\'ebrique r\'eelle de degr\'e $m$ dans l'espace projectif
complexe ${\C}P^n$ de dimension $n$ est le lieu des z\'eros dans ${\C}P^n$
d'un polyn\^ome ${\bar{f}}$ homog\`ene de degr\'e $m$ en $n+1$ variables et \`a coefficient
r\'eel (ici et par la suite ${\bar{f}}$ d\'esignera l'homog\'einis\'e d'un polyn\^ome affine $f$).
Sa partie r\'eelle est alors le lieu des z\'eros dans ${\R}P^n$ du polyn\^ome ${\bar{f}}$.
En g\'en\'eral, on notera $X_m^n$ une hypersurface alg\'ebrique r\'eelle non singuli\`ere
de degr\'e $m$ dans l'espace projectif complexe ${\C}P^n$.
Les r\'esultats pr\'esent\'es ici s'inscrivent
dans le cadre du probl\`eme de la classification des types topologiques possibles de ${\R}X_m^n$
pour $m$ et $n$ fix\'es ($16^{eme}$ probl\`eme d'Hilbert).

On appelle {\it hypersurface doubl\'ee}
toute hypersurface projective r\'eelle $X$ de degr\'e $2k$ dans ${\C}P^n$ obtenue comme
petite d\'eformation d'une hypersurface double:
$$X=\{{\bar{f}}_k(Z)^2- \epsilon \cdot {\bar{f}}_{2k}(Z)=0\} \subset {\C}P^n,$$
o\`u $\epsilon >0$ est suffisamment petit, et ${\bar{f}}_k$ (resp. ${\bar{f}}_{2k}$) est un polyn\^ome
r\'eel homog\`ene de degr\'e $k$ (resp. $2k$) en $n+1$ variables $Z_0,\cdots,Z_n$.
Dans ce papier, les hypersurfaces ${\R}X_k^n$ et ${\R}X_{2k}^n$ d\'efinies respectivement
par ${\bar{f}}_k$ et ${\bar{f}}_{2k}$ seront non singuli\`eres et s'intersecteront transversalement.
Dans ce cas ${\R}X$ est non singuli\`ere et s'obtient en ``doublant''
${\R}X_{k,+}^n={\R}X_k^n \cap \{{\bar{f}}_{2k} \geq 0\}$. On peut le voir directement ou bien
remarquer que ${\R}X$ s'obtient par une petite d\'eformation r\'eelle \`a partir de la sous
vari\'et\'e lisse  ${\R}Y=\{U^2-{\bar{f}}_{2k}(Z)=0,{\bar{f}}_k(Z)=0\}$
de l'espace projectif tordu ${\R}P^{n+1}(1,k)$ (les $Z_i$ ont poids $1$ et $U$ le poids $k$)
en consid\'erant la famille $\{U^2-{\bar{f}}_{2k}(Z)=0,{\bar{f}}_k(Z)=t \cdot U\}$ param\'etr\'ee
par $0 \leq t \leq \sqrt{\epsilon}$ (en particulier ${\R}X$ est hom\'eomorphe \`a ${\R}Y$).

Les hypersurfaces doubl\'ees ont d\'eja \'et\'e utilis\'ees
en topologie des vari\'et\'es alg\'ebriques r\'eelles. \`A titre d'exemple,
tous les types topologiques possibles de quartiques ${\R}X_4^3$ dans ${\R}P^3$
peuvent, \`a une exception pr\`es,  \^etre r\'ealis\'es comme parties r\'eelles de
quadriques doubl\'ees \cite{Vimult}. 
Les r\'esultats et constructions pr\'esent\'es ici peuvent \^etre vus comme des g\'en\'eralisations
de ceux de \cite{Basymp} qui concernaient les surfaces alg\'ebriques r\'eelles dans ${\C}P^3$
(voir la remarque \ref{gen}). On construit ici des hypersurfaces doubl\'ees de tout degr\'e et
toute dimension en utilisant la m\'ethode bien connue due \`a O. Viro
\cite{Vmet1, Vmet2, Vmet3, Ris} de
construction de vari\'et\'es alg\'ebriques r\'eelles avec topologie prescrite.
Les polyn\^omes affines $f_k$ et $f_{2k}$ d\'efinissant l'hypersurface doubl\'ee $X$
sont des polyn\^omes de Viro associ\'ees \`a des fonctions
convexes globalement affines, et l'extension de la m\'ethode de Viro au cas des intersections compl\`etes
obtenue dans \cite{Bint} permet alors de d\'eterminer la topologie du triplet $({\R}P^n,
{\R}X_k^n, {\R}X_{2k}^n)$. On obtient que
l'hypersurface doubl\'ee ${\R}X$ est hom\'eomorphe au r\'esultat d'un ``d\'ecoupage-collage''
sur l'union des vari\'et\'es produits ${\R}X_k^l \times {\R}Y_{2k}^{n-l}$, $l=1,\cdots,n$, o\`u
$X_k^l$ est l'intersection
de $X_k^n$ avec un $l$-sous-espace de coordonn\'ees de ${\C}P^n$ et $Y_{2k}^{n-l}$
est le rev\^etement double r\'eel d'un $(n-l)$-sous-espace de coordonn\'ees
ramifi\'e sur son intersection avec $X_{2k}^n$.
De tels rev\^etements doubles $Y_{2k}^{n-l}$ sont appel\'es
{\it plans doubl\'es} et occupent une place importante
dans l'\'etude des vari\'et\'es r\'eelles, notamment (pour $n-l=2$)
des courbes r\'eelles planes. Notons qu'un plan doubl\'e r\'eel
$Y_{2k}^{n}$ se plonge dans l'espace projectif tordu ${\C}P^{n+1}(1,k)$
(muni de la conjugaison complexe usuelle)
comme $Y_{2k}^{n}=\{U^2 \pm {\bar{f}}_{2k}(Z)=0\}$ si $X_{2k}^n=\{{\bar{f}}_{2k}(Z)=0\}$
est le lieu de ramification (on dira que $\{U^2-{\bar{f}}_{2k}(Z)=0\}$
est le plan doubl\'e r\'eel associ\'e au polyn\^ome $f_{2k}$).
Il est en g\'en\'eral assez difficile (pour $n \geq 3$) d'expliciter compl\^etement
la topologie de ${\R}X$ du fait des collages. N\'eanmoins, on peut se concentrer
sur la ``partie asymptotique'' du type topologique de ${\R}X$ lorsque $k \rightarrow +\infty$,
auquel cas les collages deviennent n\'egligeables et
les hypersurfaces $X_k^l$ et $X_{2k}^{n-l}$ peuvent \^etre
choisies arbitrairement (Proposition \ref{cle}). On peut alors obtenir
des estimations sur les comportements asymptotiques pour $m$ (resp. $k$) $\rightarrow +\infty$
des valeurs maximales des nombres de betti $b_i({\R}X_m^n)$  (resp. $b_i({\R}Y_{2k}^n)$)
lorsque $i$ et $n$ sont fix\'es (les groupes d'homologie sont pris \`a coefficients
dans ${\Z}/2$). On montre (Proposition \ref{as1}) que les suites $Max\,  b_i({\R}X_m^n)$ et $Max \,  b_i({\R}Y_{2k}^n)$
index\'ees par $m$ et $k$, sont asymptotiquement \'equivalentes \`a ${\h}_{i,n} \cdot m^n$
et \`a ${\d}_{i,n} \cdot k^n$, respectivement, pour certains r\'eels ${\h}_{i,n}$ et ${\d}_{i,n}$.
En d'autres termes, pour $i$ et $n$ fix\'es, on montre l'existence des limites

$$\lim_{m \rightarrow +\infty} \frac{Max\,  b_i({\R}X_m^n)}{m^n}={\h}_{i,n}, \quad
\lim_{k \rightarrow +\infty} \frac{Max\,  b_i({\R}Y_{2k}^n)}{k^n}={\d}_{i,n}.$$

Notre construction d'hypersurfaces doubl\'ees permet alors
de montrer les in\'egalit\'es (th\'eor\`emes \ref{mth} et et \ref{mth1})

$$ {\h}_{0,n} \geq \frac{1}{2^n-2}
\sum_{l=1}^{n-1} {\h}_{0,l} \cdot {\d}_{0,n-l},
\quad {\d}_{0,n} \geq  \sum_{l=1}^{n-1} {\h}_{0,l} \cdot {\d}_{0,n-l} \; + {\h}_{0,n}.$$

On obtient ainsi des bornes inf\'erieures r\'ecursives pour les ${\h}_{0,n}$ et ${\d}_{0,n}$, desquelles est extraite une borne
explicite (Proposition \ref{le})

$${\h}_{0,n} \geq \frac{1}{2^{n-1}}.$$

Des bornes sup\'erieures (Proposition \ref{bsup}) sont \'egalement obtenues de mani\`ere classique
en utilisant l'in\'egalit\'e de Smith-Thom
$$b_*({\R}X) \leq b_*(X),$$ valable pour toute vari\'et\'e alg\'ebrique r\'eelle $X$ et o\`u $b_*$
d\'esigne la somme totale des nombres de Betti, ainsi que les in\'egalit\'es de Comessatti-Petrowsky-Oleinik
g\'en\'eralis\'ees
$$|\chi({\R}X)-1| \leq h^{\frac{n}{2},\frac{n}{2}}(X)-1,$$
valables (en particulier) pour toute vari\'et\'e alg\'ebrique projective r\'eelle non singuli\`ere de dimension paire $n$ et o\`u
$\chi$ d\'esigne la caract\'eristique d'Euler et $ h^{\frac{n}{2},\frac{n}{2}}$ le nombre de Hodge de type $(\frac{n}{2},\frac{n}{2})$.
On renvoie \`a \cite{Kha} pour un survey r\'ecent sur la topologie
des vari\'et\'es alg\'ebriques r\'eelles.

Le tableau 1. (sous-section \ref{est}) recense les estimations obtenues pour les
${\h}_{0,n}$ et ${\d}_{0,n}$ tels que $n \leq 7$.

On s'int\'eresse ensuite au cas des surfaces alg\'ebriques r\'eelles dans ${\C}P^3$.
On montre les in\'egalit\'es (Th\'eor\`eme \ref{surfg})

$$ \frac{{\d}_{0,2}}{6}+\frac{1}{12} \leq {\h}_{0,3} \leq \frac{5}{12}, \quad \frac{{\d}_{1,2}}{6}+\frac{1}{6}
\leq {\h}_{1,3} \leq \frac{5}{6}.$$

Comme corollaire de ces in\'egalit\'es et des in\'egalit\'es ${\d}_{0,2} \geq \frac{27}{16}$ et ${\d}_{1,2}\geq \frac{27}{8}$
r\'ealis\'ees par de r\'ecents contre-exemples \`a la conjecture
de Ragsdale obtenus par Itenberg \cite{Itrag}, on am\'eliore (Th\'eor\`eme \ref{surf})
le r\'esultat principal de \cite{Basymp} (voir la remarque \ref{gen})

$$ \frac{35}{96} \leq {\h}_{0,3} \leq \frac{5}{12}, \quad \frac{35}{48}
\leq {\h}_{1,3} \leq \frac{5}{6}.$$

On s'int\'eresse ensuite (voir la sous-section \ref{T-hyp})
aux $T$-vari\'et\'es alg\'ebriques r\'eelles qui sont par d\'efinition
des vari\'et\'et\'es alg\'ebriques r\'eelles construites par la version combinatoire de la m\'ethode de Viro
appel\'ee {\it patchwork combinatoire} ou {\it T-construction} (voir la sous-section \ref{T-hyp}).
Le patchwork combinatoire est un outil de construction tr\`es puissant. Il a par exemple servi \`a
Itenberg pour obtenir des contre-exemples \`a une vieille conjecture attribu\'ee \`a Ragsdale (voir \cite{ItT-cour,Itrag}),
ainsi qu'\`a Viro et Itenberg \cite{ItVi} pour montrer l'exactitude de l'in\'egalit\'e de Smith-Thom pour les hypersurfaces alg\'ebriques
$X_m^n \subset {\C}P^n$ de tout degr\'e et toute dimension. La question est alors d'estimer la ``richesse'' des $T$-vari\'et\'es.
On montre (Proposition \ref{as1T}) l'existence de nombres r\'eels $Th_{i,n}$ et $Td_{i,n}$ d\'efinis comme pour
$h_{i,n}$ et $d_{i,n}$ mais en se restreignant aux $T$-vari\'et\'es de la famille correspondante. Des bornes sup\'erieures
pour certains de ces  nombres ont r\'ecemment \'et\'e obtenues par Itenberg et Shustin \cite{ItSh}. Notre construction permet alors de montrer
(Th\'eor\`eme \ref{Tsurf}) l'existence d'hypersurfaces $X_m^n \subset {\C}P^n$ qui ne sont pas des $T$-hypersurfaces
pour tout $n \geq 5$ et tout degr\'e $m$ suffisamment grand.
\smallskip

Le papier est divis\'e en trois sections. La premi\`ere section est consacr\'ee
\`a des rappels sur la m\'ethode de viro.
Dans la deuxi\`eme section on pr\'esente notre construction d'hypersurfaces doubl\'ees.
La derni\`ere section est consacr\'ee aux applications de cette construction concernant des
valeurs maximales asymptotiques de nombres de Betti.

\section{M\'ethode
de Viro}

On commence par rappeler la notion de carte d'un polyn\^ome.

\subsection{Cartes
d'un polyn\^ome}
\label{carte}

\`A partir de maintenant, un point entier de ${\R}^n$ est
un point \`a coordonn\'ees enti\`eres et un polytope est
un polytope convexe \`a sommets entiers dans ${\R}^n$.
Soit $T_m^n$ le polytope de Newton d'un polyn\^ome affine
g\'en\'erique de degr\'e $m$ en $n$ variables:
$$T_m^n=conv \, \{(y_1,\cdots,y_n) \in {\R}^n, \; y_1,\cdots,y_n \geq 0,
\; y_1+\cdots+y_n \leq m\}.$$
On note $e_0=(0,\cdots,0)$ et $e_i$ le sommet de $T_m^n$
dont la coordonn\'ee $y_i$ est \'egale \`a $m$.
Soient $Z=(Z_0:\cdots:Z_n)$ des coordonn\'ees homog\`enes
sur ${\C}P^n$. On fait correspondre \`a chacune des faces $\Gamma$ de $T_m^n$ (comprenant $T_m^n$
lui m\^eme) un sous-espace projectif de ${\C}P^n$
comme suit: si $\Gamma=conv\{e_i, \, i \in I\}$ avec $I \subset \{0,\cdots,n\}$,
alors $X(\Gamma)$ est le sous-espace projectif ${\C}P^{|I|-1}$
de coordonn\'ees homog\`enes $Z_i, \, i \in I$, i.e. le sous-espace d'\'equations $Z_j=0, \; j \in \{0,\cdots,n\} \setminus I$.

On consid\`ere la carte affine $\{Z_0 \neq 0\}$ munie des coordonn\'ees
affines $z=(z_1,\cdots,z_n)$ avec $z_i=Z_i/Z_0, \, i=1,\cdots,n$.
Si $f(z)$ est un polyn\^ome de polytope de newton $T_m^n$ et $\Gamma$ est une face de $T_m^n$, alors
l'intersection de l'hypersurface $\{{\bar{f}}=0\} \subset {\C}P^n$
(d\'efinie par l'homog\'einis\'e de $f$) et du sous-espace $X(\Gamma)$ coincide
avec l'hypersurface $\{{\bar{f}}^{\Gamma}=0\} \subset X(\Gamma)$, o\`u $f^{\Gamma}$ est le
tronqu\'e de $f$ sur $\Gamma$: $f^{\Gamma}(z) =\sum_{w \in \Gamma} a_wz^w$ si $f(z)=\sum a_wz^w$.

Ce que l'on vient de d\'ecrire est un cas particulier
d'une correspondance $P \rightarrow X(P)$ (surjective mais non injective)
de l'ensemble des polytopes sur celui des vari\'et\'es toriques projectives
(normales ou non) ayant les propri\'et\'es suivantes.
Pour toute face $\Gamma$ d'un polytope $P$ la vari\'et\'e
$X(\Gamma)$ se plonge comme une sous vari\'et\'e torique de $X(P)$ de telle
sorte que pour deux faces $\Gamma$ et $\Gamma'$ de $P$ l'on ait
$X(\Gamma) \cap X(\Gamma')= X(\Gamma \cap \Gamma')$.
Tout polyn\^ome $f$ de polytope de Newton $P$ d\'efinit une hypersurface
de $X(P)$, cette hypersurface intersecte $X(\Gamma)$ le long de l'hypersurface d\'efinie par
$f^{\Gamma}$ pour toute face $\Gamma$ de $P$.

Soit $P$ un polytope situ\'e dans l'orthant positif ${({\R}_+)}^n$.
L'application moment associ\'ee \`a un ensemble de points entiers $\mathcal P$
d'enveloppe convexe $P$ est l'application $\phi:{({\C}^*)}^n \rightarrow P$
d\'efinie par
$$\phi(z)=\frac{\sum_{i=0}^s {|z^{w_i}|w_i}}{\sum_{i=0}^s {|z^{w_i}|}}.$$

On identifie \`a ${({\Z}/2)}^n$ le groupe des sym\'etries de ${\R}^n$ par rapport aux hyperplans de coordonn\'ees
$\{y_i=0\}$ ainsi que le groupe des sym\'etries de ${\R}^n$ par rapport aux hyperplans de coordonn\'ees
$\{z_i=0\}$ via les isomorphismes $g=(\epsilon_1,\cdots,\epsilon_n) \mapsto \big((y_1,\cdots,y_n) \mapsto ((-1)^{\epsilon_1}y_1,
\cdots,(-1)^{\epsilon_n}y_n)\big)$ et $g=(\epsilon_1,\cdots,\epsilon_n) \mapsto \big((z_1,\cdots,z_n) \mapsto ((-1)^{\epsilon_1}z_1,
\cdots,(-1)^{\epsilon_n}z_n)\big)$. Le groupe de sym\'etries identifi\'e \`a ${({\Z}/2)}^n$ sera clair suivant le contexte.

On note ${\R}(g)$ l'orthant obtenu en prenant l'image de ${({\R}^*_+})^n$
par $g \in {({\Z}/2)}^n$. On consid\`ere la restriction de $\phi$ \`a ${({\R}^*_+})^n$ et
on \'etend l'application obtenue en une application (que l'on appellera
\`a nouveau application moment) $\tilde{\phi}:{({\R}^*)}^n \rightarrow P^*$ o\`u
$P^*=\cup_{g \in
{({\Z}/2)}^n}g(P)$ par la r\`egle $\tilde{\phi}(g(z))
=g(\phi(z))$ si $z \in {({\R}^*_+})^n$.
Si $P$ est d'int\'erieur non vide i.e. $\dim(P)=n$,
alors $\tilde{\phi}$ est un diff\'eomorphisme de
${\R}(g)$ sur l'int\'erieur de $g(P)$ pour tout
$g \in {({\Z}/2)}^n$.
Soit $C_P$ la vari\'et\'e topologique r\'esultant des identifications
suivantes sur $P^*$:
pour toute face $\Gamma$ de $P^*$ et tout vecteur
$\alpha \in {\Z}^n$ orthogonal \`a $\Gamma$,
la face $\Gamma$ est identifi\'ee \`a la copie sym\'etrique $g(\Gamma)$
o\`u $g$ est la r\'eduction de $\alpha$ dans ${({\Z}/2)}^n$.
Pour toute face $\Gamma$ de $P$, l'image de $\cup _{g \in {({\Z}/2)}^n}
g(\Gamma)$ dans $C_P$ coincide avec $C_{\Gamma}$
et, si $\dim(P)=n$, l'application $\tilde{\phi}$ peut \^etre \'etendue en
un hom\'eomorphisme ${\R}X(P) \rightarrow C_P$ envoyant ${\R}X(\Gamma)$
sur $C_{\Gamma}$ pour toute face $\Gamma$ de $P$. Un tel hom\'eomorphisme
sera dit stratifi\'e par la suite.

A titre d'exemple, si $P=T_m^n$ alors $C_P$ s'obtient
en identifiant les points antipodaux situ\'es sur le bord du polytope
$$conv \,  \{(\pm m,0,\cdots,0),(0,\pm m,0,\cdots,0),\cdots,(0,\cdots,0,\pm m)\}.$$

Soit $f$ un polyn\^ome de polytope de Newton $P$. On appelle carte de $f$
chacun des ensembles
$$C_P^{g}(f) :=\tilde{\phi}(\{z \in {\R}(g), \;
f(z)=0\}) \subset g(P), \quad g \in {({\Z}/2)}^n,$$
ainsi que l'adh\'erence $C_P(f)$ de
$\tilde{\phi}(\{z \in {({\R}^*)^n}, \; f(z)=0\})$ dans $C_P$
(ces cartes d\'ependent de l'ensemble $\mathcal P$ choisi, dans ce papier
on prendra toujours ${\mathcal P}=P \cap {\Z}^n$).
Notons que l'on a $C_P^{g}(f)=g(C_P^{id}(f \circ g))$.
La carte $C_P(f)$ s'obtient par recollement des cartes
$C_P^{g}(f)$ dans $C_P$ et v\'erifie la propri\'et\'e
$C_P(f) \cap C_{\Gamma}=C_{\Gamma}(f^{\Gamma})$ pour toute face $\Gamma$ de $P$.
Si $\dim(P)=n$, l'application $\tilde{\phi}$ donne, pour tout $g
\in {({\Z}/2)}^n$, un hom\'eomorphisme de paire
$({\R}(g),\{f=0\} \cap {\R}(g)) \simeq (int(g(P)), C_P^{g}(f))$
et s'\'etend en un hom\'eomorphisme $({\R}X(P),\{{\bar{f}}=0\}) \simeq (C_P, C_P(f))$.
Le polyn\^ome $f$ est dit {\it non d\'eg\'en\'er\'e} si pour toute face $\Gamma$ de $P$
(comprenant $P$ lui m\^eme) le tronqu\'e de $f$ sur $\Gamma$
d\'efinit une hypersurface non singuli\`ere de ${({\R}^*)^n}$. Si $f$
est non d\'eg\'en\'er\'e et $X(P)$ est non singuli\`ere
alors l'hypersurface de ${\R}X(P)$ d\'efinie par $f$ est non singuli\`ere.
\smallskip

\subsection{Transformations affines enti\`eres unimodulaires}
\label{Trans}

On rappelle des r\'esultats bien connus qui nous seront utiles \`a plusieurs reprises.
Notons $AFF_n(\Z)$ le groupe des transformations affines unimodulaires de ${\R}^n$
\`a coefficients entiers. Soit $\Delta \in AFF_n(\Z)$. La transformation $\Delta$ est la compos\'ee
d'une translation par un point $b \in {\Z}^n$ et d'une transformation $L \in GL_n(\Z)$,
que l'on appellera partie lin\'eaire de $\Delta$.
Le changement de
coordonn\'ees multiplicatif du tore complexe correspondant \`a $\Delta$ (en fait \`a $L$) est l'application
$$\Delta_*:{({\C}^*)}^n \rightarrow {({\C}^*)}^n, \quad
z=(z_1,\cdots,z_n) \mapsto \tilde{z}=(\tilde{z}_1, \cdots ,\tilde{z}_n)$$
d\'efinie par
$$\tilde{z_j}=\prod_{i=1}^n z_i^{{(L^{-1})}_{ij}}, \quad j=1,\cdots,n.$$

Soit $P$ un polytope dans ${\R}^n$. On a alors $X(\Delta(P))=X(P)$.
L'image par $\Delta$ d'un polyn\^ome $f \in \R[z_1,\cdots,z_n]$ d\'efini par $f(z)=\sum a_w z^w$ est le polyn\^ome
$\tilde{f} \in \R[\tilde{z}_1, \cdots ,\tilde{z}_n]$ d\'efini par
$\tilde{f}(\tilde{z})={\tilde{z}}^b \cdot \sum a_w {\tilde{z}}^{L(w)}$.
Si $P$ est le polytope de Newton de $f$, alors le polytope de Newton de $\tilde{f}$ est $\Delta(P)$.
\`A partir de l'\'egalit\'e $z^w={\tilde{z}}^{L(w)}$ valable pour tout $w \in {\Z}^n$, on obtient facilement que
$f$ et $\tilde{f}$ d\'efinissent la m\^eme hypersurface
dans ${({\C}^*)}^n$, et donc dans $X(\Delta(P))=X(P)$. En particulier, $f$ est non d\'eg\'en\'er\'e
si et seulement si $\tilde{f}$ l'est. De plus, si $\phi:{({\C}^*)}^n \rightarrow P$ est l'application moment
relative au syst\`eme de coordonn\'ees $(z_1,\cdots,z_n)$ et ${\tilde{\phi}}:{({\C}^*)}^n \rightarrow \Delta(P)$
est l'application moment relative au syst\`eme de coordonn\'ees $(\tilde{z}_1, \cdots ,\tilde{z}_n)$, alors on
a $(\Delta \circ \phi)(z)=\tilde{\phi}(\tilde{z})$. En utilisant cette \'egalit\'e, on peut toujours se ramener
au cas o\`u le polytope est non vide (de dimension maximale).

\subsection{M\'ethode de Viro pour les hypersurfaces}
\label{hypvir}

On consid\`ere une hypersurface r\'eelle de degr\'e $m$ dans ${\C}P^n$
d\'efinie, pour $t>0$ suffisamment petit, par un polyn\^ome affine
$$
f_t(z)=\sum_{w \in {\mathcal T}}a_{w}t^{\nu(w)}z^w
$$
de polytope de Newton $T_m^n$, o\`u les $a_{w}$ sont des coefficients r\'eels et $\nu$
est une application \`a valeurs enti\`eres d\'efinie sur ${\mathcal T} \subset {(\Z)}^n$.
Un tel polyn\^ome est appel\'e {\it polyn\^ome de Viro}.

On associe \`a $\nu$ une subdivision poly\'edrale $\tau$ de $T$
en projettant par $\pi:{\R}^n \times {\R} \rightarrow {\R}$ la r\'eunion
$G$ des faces de la partie inf\'erieure du bord du polytope 
$$\hat{T}=conv \,  \{(w,\nu(w)), \; w \in {\mathcal T}\}$$
($G$ est l'union des faces compactes du poly\`edre $(0,\cdots,0) \times {\R}_++\hat{T}$).
Une subdivision poly\'edrale obtenue de cette mani\`ere est dite {\it convexe} (ou coh\'erente).
Soit $f$ le polyn\^ome
obtenu en posant $t=1$ dans $f_t$.
On fait l'hypoth\`ese suivante.
\medskip

{\bf Hypoth\`ese.} Pour tout $F \in \tau$
le polyn\^ome $f^F$ est non d\'eg\'en\'er\'e.
\medskip

Les cartes $C_{F}^{g}(f^F)$,
$F \in \tau$ et $g \in {({\Z}/2)}^n$, se recollent entre elles
dans $C_T$. Soit $S$ le r\'esultat d'un tel collage.
Notons $X$ l'hypersurface r\'eelle $\{{\bar{f}}_t(x)=0\}$ de degr\'e $m$ dans 
${\C}P^n$.

Rappelons qu'un hom\' eomorphisme ${\R}X(P) \rightarrow C_{P}$
est dit stratifi\'e si il envoie ${\R}X(\Gamma)$ sur $C_{\Gamma}$
pour toute face $\Gamma$ de $P$. On peut maintenant \'enoncer le th\'eor\`eme de Viro
\cite{Vmet1,Vmet2,Vmet3} (voir aussi \cite{Ris, Bint}).

\begin{thm}[Viro]
\label{methyp}
Pour $t>0$ suffisamment petit le polyn\^ome $f_t$
est non d\'eg\'en\'er\'e (en particulier ${\R}X$ est non singuli\`ere)
et il existe un hom\'eomorphisme stratifi\'e
${\R}P^n \rightarrow C_{T}$ envoyant ${\R}X$
sur $S$.
\end{thm}

\subsection{M\'ethode de Viro pour les hypersurfaces doubl\'ees}
\label{doubcons}

On consid\`ere une hypersurface doubl\'ee r\'eelle
$X=\{{{\bar f}_k}^2-\epsilon \cdot {\bar f}_{2k}=0\}$ ($0< \epsilon \ll 1$)
de degr\'e $2k$ dans ${\C}P^n$ o\`u $f_k$ et $f_{2k}$ sont, pour $t>0$ suffisamment petit,
des polyn\^omes de Viro $f_k=f_{1,t}$ et $f_{2k}=f_{2,t}$
$$f_{i,t}(z)=\sum_{w \in {\mathcal T}_i}a_{i,w}t^{\nu_i(w)}z^w, \quad i=1,2,$$
de polytopes de Newton $T_1=T_{k}^n$ et $T_2=T_{2k}^n$, respectivement.

On note $\tau_i$ la subdivision poly\'edrale de $T_i$
associ\'ee \`a $\nu_i$. Soit $T=T_{3k}^n$
la somme de Minkowsky $T_1+T_2$. On associe \`a la paire
$(\nu_1,\nu_2)$ une subdivision poly\'edrale $\tau$ de $T$
en projettant par $\pi:{\R}^n \times {\R}\rightarrow {\R}^n$
la r\'eunion $G$ des faces inf\'erieures du polytope
$$\hat{T}=conv \, \{(w_1+w_2,\nu_1(w_1)+\nu_2(w_2)), \; w_i \in {\mathcal T}_i\}.$$
Notons $G_i$ la r\'eunion des faces de la partie inf\'erieure
du polytope $\hat{T}_i$ correspondant \`a $\nu_i$ ($\tau_i$ s'obtient donc
en projettant $G_i$). Alors chacune des facettes (faces de dimension maximale)
de $G$ s'\'ecrit de mani\`ere unique comme la somme de Minkowsky d'une face de $G_1$ et
d'une face de $G_2$. En projettant, on obtient alors une repr\'esentation, induite par $(\nu_1,\nu_2)$,
de chaque polytope $F \in \tau$ comme $F=F_1+F_2, \quad F_i \in \tau_i$.
Par la suite, lorsque l'on \'ecrira $F=F_1+F_2$ avec $F \in \tau$ et $F_i \in \tau_i$,
on se r\'eferrera toujours \`a la repr\'esentation induite par $(\nu_1,\nu_2)$.
Pour $i=1,2$, on note $f_i$ le polyn\^ome obtenu en posant $t=1$
dans $f_{i,t}$. On fait les hypoth\`eses suivantes.
\medskip

{\bf Hypoth\`eses.}
\begin{enumerate}
\item
Pour tout $F_i \in \tau_i$,
le polyn\^ome $f_i^{F_i}$ est non d\'eg\'en\'er\'e.
\item  La paire $(\nu_1,\nu_2)$ est
suffisamment g\'en\'erique au sens o\`u si $F=F_1+F_2$ avec $F \in \tau$ et $F_i \in \tau_i$,
alors $\dim(F)=\dim(F_1)+\dim(F_2)$.
\end{enumerate}
Sous la deuxi\`eme hypoth\`ese, la subdivision $\tau$ est dite {\it mixte}.
\medskip

Pour chaque polytope $F \in \tau$ et chaque $g \in {({\Z}/2)}^n$,
si $F=F_1+F_2$ avec $F_i \in \tau_i$,
on consid\`ere les sous-ensembles de $g(F)= g(F_1)+g(F_2)$
$$C_{F_1}^{g}(f_1^{F_1})+g(F_2), \quad
g(F_1)+C_{F_2}^{g}(f_2^{F_2}).$$
Les ensembles $C_{F_1}^{g}(f_1^{F_1})+g(F_2)$, $F \in \tau$ et $g \in {({\Z}/2)}^n$, se recollent 
entre eux dans $C_T$. Soit $S^1$ le r\'esultat d'un tel collage.
De m\^eme, les ensembles $g(F_1)+C_{F_2}^{g}(f_2^{F_2})$, $F \in \tau$ et $g \in {({\Z}/2)}^n$, se recollent
entre eux dans $C_T$. Soit $S^2$ le r\'esultat d'un tel collage. Notons
$X_k$ et $X_{2k}$ respectivement les hypersurfaces $\{{\bar f}_{1,t}=0\}$
et $\{{\bar f}_{2,t}=0\}$ de ${\C}P^n$.

Le r\'esultat suivant est un corollaire imm\'ediat de (\cite{Bint}, Th\'eor\`eme 2.1.)
\begin{prop}
\label{int}
Pour $t>0$ suffisamment petit, les hypersurfaces ${\R}X_k$ et ${\R}X_{2k}$
sont non singuli\`eres, s'intersectent transversalement, et
il existe un hom\'eo\-morphisme stratifi\'e
${\R}P^n \rightarrow C_{T}$ envoyant ${\R}X_k$
sur $S^1$ et ${\R}X_{2k}$ sur $S^2$.
\end{prop}

Notons que la proposition pr\'ec\'edente implique que, pour $t>0$ et $\epsilon>0$ suffisamment petits,
l'hypersurface doubl\'ee ${\R}X=\{{\bar f}_{1,t}^2-\epsilon \cdot {\bar f}_{2,t}=0\} \subset {\R}P^n$ est non singuli\`ere.
On d\'efinit deux sous-espaces $S^2_+$ et $S^2_-$ de $C_T$ par
les conditions
$$C_{T}=S^2_+ \cup S^2_-, \quad S^2=S^2_+ \cap S^2_-$$
et, si $\mathbf 0$ d\'esigne l'image dans $C_T$ du sommet $(0,\cdots,0)$ de $T$,
$${\mathbf 0} \in S^2_{\pm} \iff \pm f_{2,t}(0,\cdots,0) \geq 0.$$

Un hom\'eomorphisme stratifi\'e ${\R}P^n \rightarrow C_{T}$
envoie $(1:0:\cdot:0)$, qui est la vari\'et\'e torique associ\'ee au sommet $(0,\cdots,0)$ de $T$,
sur $\mathbf 0$. Par suite, un tel hom\'eomorphisme envoie $\{\pm {\bar f}_{2,t}(Z) \geq 0\}$ sur $S^2_{\pm}$.
Rappelons que ${\R}X$ s'obtient en doublant ${\R}X_{k,+}:=\{Z \in {\R}X_k, \;
{\bar f}_{2,t}(Z) \geq 0\}$. La proposition \ref{int} implique donc le r\'esultat suivant.
\begin{cor}
\label{hom}
Pour $t, \epsilon>0$ suffisamment petits, l'hypersurface doubl\'ee ${\R}X$ est non singuli\`ere
et ${\R}X_{k,+}$ est hom\'eomorphe \`a $S^2_+ \cap S^1$.
\end{cor}

Soit $C_{F_2,+}^{g}(f_2^{F_2}) \subset g(F_2)$
l'image de $\{z \in {\R}(g), \; f_2^{F_2}(z) \geq 0\}$ par
une application moment $\tilde{\phi}$ envoyant $\{z \in {\R}(g), \; f_2^{F_2}(z)=0\}$
sur $C_{F_2}^{g}(f_2^{F_2})$.
Alors $S^2_+$ s'obtient par recollement des
$$g(F_1)+C_{F_2,+}^{g}(f_2^{F_2}) \subset g(F)$$
et donc $S^2_+ \cap S^1$ s'obtient par recollement des
$$C_{F_1}^{g}(f_1^{F_1})+C_{F_2,+}^{g}(f_2^{F_2}) \subset g(F).$$

\section{Construction}
\label{cons}

On reprend la construction d\'ecrite dans la sous-section
\ref{doubcons} d'une hypersurface doubl\'ee.
\smallskip

Soit  $X=\{{{\bar f}_k}^{\, 2}-\epsilon \cdot {\bar f}_{2k}=0\}$ ($0< \epsilon \ll 1$)
une hypersurface doubl\'ee de degr\'e $2k$ dans ${\C}P^n$ o\`u $f_k$ et $f_{2k}$ sont donn\'es, pour $t>0$ suffisamment petit,
par des polyn\^omes de Viro $f_k=f_{1,t}$ et $f_{2k}=f_{2,t}$
$$f_{i,t}(z)=\sum_{w \in {\mathcal T}_i}a_{i,w}t^{\nu_i(w)}z^w, \quad i=1,2$$
et supposons \`a partir de maintenant que, pour $i=1$ et $i=2$,
\smallskip

\begin{center}
{\bf l'application $\nu_i$ est la restriction d'une application affine sur ${\R}^n$}.
\end{center}
\smallskip

Chacune des triangulations $\tau_i$ associ\'ee \`a $\nu_i$ est alors triviale
au sens o\`u elle consiste en la r\'eunion des faces de $T_i$ (en particulier, la triangulation
$\tau_i$ contient $T_i$ comme seul $n$-simplexe et les sommets de $\tau_i$ sont les sommets de $T_i$).
Pour $i=1,\cdots,n$, on note $u_i$ (resp. $v_i$) le sommet de $T_1$ (resp. $T_2$) dont la i-\`eme
coordonn\'ee est \'egale \`a $k$ (resp. $2k$) puis $u_0=v_0=(0,\cdots,0)$.
Soit $\tau$ la subdivision convexe de $T=T_1+T_2$ associ\'ee \`a la paire
$(\nu_1,\nu_2)$. Chacun des polytopes de $\tau$ est la forme $F_1+F_2$
o\`u $F_i$ est une face de $T_i$. Rappelons que $\tau$ est mixte si et seulement si
chacune de ces sommes $F_1+F_2$ est une somme directe.

\begin{lem}
\label{cayley trick}
La subdivision $\tau$ est une subdivision mixte 
si et seulement s'il existe
une permutation $\sigma$ de $\{0,\cdots,n\}$ pour laquelle
$\tau$ soit consitu\'ee des $n$-polytopes
$$F_l=conv\{u_{\sigma(0)},\cdots,u_{\sigma(l)}\}+conv\{v_{\sigma(l)},\cdots,v_{\sigma(n)}\}, \quad l=0,\cdots, n$$
(et de leurs faces).
\end{lem}

{\sc Preuve.}
Le ``Cayley trick'' combinatoire (voir \cite{Stmix,Hmix})
\'etablit une bijection entre les subdivisions mixtes convexes de $T_1+T_2$ \`a sommets dans $\{u_i+v_j \, | \, i,j=0,\cdots,n\}$
et les triangulations convexes \`a sommets dans
$\{u_i \times \{0\} \, | \, i=0,\cdots,n\} \cup \{v_j \times \{1\} \, | \, i=0,\cdots,n\}$ du prisme obtenu en prenant le join
$\big(T_1 \times \{0\}\big) \star \big(T_2 \times \{1\}\big) \subset {\R}^{n+1} $. Cette bijection envoie un simplexe
$\big(\sigma_1 \times \{0\}\big)  \star \big(\sigma_2 \times \{2\}\big) $
de la triangulation du prisme ($\sigma_i$ est une face de $T_i$) sur le polytope
$\sigma_1+\sigma_2$. De mani\`ere \'evidente, l'ensemble des triangulations convexes du
type pr\'ec\'edent est en bijection avec celui des triangulations convexes \`a sommets entiers du prisme $T_1^n \times T_1^1$.
Le r\'esultat d\'ecoule alors de la description connue des triangulations convexes
\`a sommets entiers de $T_1^n \times T_1^1$ (voir, par exemple, \cite{GKZ}, section 7.3.C).
{\cqfd}

Pour $i=1,\cdots,n$, posons $a_i=\nu_1(u_i)$ et $b_i=\nu_2(v_i)$ de telle sorte que
$\nu_1$ et $\nu_2$ soient les restrictions des applications affines
$(y_1,\cdots,y_n) \mapsto \frac{1}{k}\sum_{i=0}^n (a_i-a_0)y_i \; +a_0$ et $(y_1,\cdots,y_n)=\frac{1}{2k}\sum_{i=0}^n
(b_i-b_0)y_i \; +b_0$, respectivement.

\begin{prop}
\label{mixte}
La subdivision $\tau$ est une subdivision mixte si et seulement si
$2a_i-b_i \neq 2a_j-b_j$ pour tout $i \neq j$, auquel cas, $\tau$ est constitu\'ee
des $n$-polytopes
$$F_l=conv\{u_{\sigma(0)},\cdots,u_{\sigma(l)}\}+conv\{v_{\sigma(l)},\cdots,v_{\sigma(n)}\}, \quad l=0,\cdots, n$$
o\`u $\sigma$ est l'unique permutation de $\{0,\cdots,n\}$
telle que
$$2a_{\sigma(0)}-b_{\sigma(0)}<2a_{\sigma(1)}-b_{\sigma(1)}<\cdots<2a_{\sigma(n)}-b_{\sigma(n)}.$$
\end{prop}

{\sc Preuve.}
Supposons que $\tau$ soit une subdivision mixte associ\'ee \`a une permutation $\sigma$ et
montrons que l'on a $2a_{\sigma(0)}-b_{\sigma(0)}<2a_{\sigma(1)}-b_{\sigma(1)}<\cdots<2a_{\sigma(n)}-b_{\sigma(n)}$.
Quitte \`a utiliser la transformation de $AFF_n(\Z)$ envoyant
$u_{\sigma(i)}+v_{\sigma(i)}$ sur $u_i+v_i$ pour $i=0,\cdots,n$
(une telle transformation envoie $u_{\sigma(i)}$ et $v_{\sigma(i)}$ sur $u_i$ et $v_i$, respectivement),
on peut se ramener au cas o\`u $\sigma$ est l'identit\'e
i.e. $\tau$ est la subdivision mixte dont les $n$-polytopes sont les
$F_l=conv\{u_0,\cdots,u_l\}+conv\{v_l,\cdots,v_n\}, \quad l=0,\cdots, n$.
Soit $\nu:T \rightarrow {\R}$ la fonction (convexe, affine par morceaux) dont le graphe est $G$.
Pour tout $l=0,\cdots,n$, il existe une fonction
affine $\alpha_l:{\R}^n \rightarrow \R$ telle que
$\nu(y) > \alpha_l(y) \;  \forall \, y \in T \setminus F_l$ et $\nu(y)= \alpha_l(y) \;  \forall \, y \in F_l$. 
Sachant que $u_i+v_j \neq u_{i'}+v_{j'} \; \forall \, (i,j) \neq (i',j')$, on obtient que $\nu(u_i+v_j)= \nu_1(u_i)+\nu_2(v_j)
=a_i+b_j$ pour tout $(i,j) \in \{0,\cdots,n\}^2$. On connait alors les valeurs de $\alpha_l$ aux sommets
de $F_l$. Ces valeurs d\'eterminent $\alpha_l$ puisque $\alpha_l$ est affine et que $F_l$ est de dimension $n$.
On obtient que $\alpha_l$ est l'application envoyant $(y_1,\cdots,y_n)$ sur
$$
\frac{1}{k}\sum_{i=1}^l(a_i-a_0)y_i+
\frac{1}{2k}\sum_{i=l+1}^n[b_i-b_l+2(a_l-a_0)]y_i
+3a_0+b_l-2a_l.$$
On montre ensuite facilement que les in\'egalit\'es
$a_i+b_j=\nu(u_i+v_j) > \alpha_l(u_i+v_j)$ valables pour tout
$l=0,\cdots,n$ et tout $ u_i+v_j  \in T \setminus F_l$ impliquent
que $2a_0-b_0<2a_1-b_1<\cdots<2a_n-b_n$.

R\'eciproquement, supposons que $\tau$ ne soit pas une subdivision mixte.
Il existe alors une face $F_1$ de $T_1$ et une face $F_2$ de $T_2$
telles que $F=F_1+F_2 \in \tau$ et $\dim(F)< \dim(F_1)+\dim(F_2)$.
Pour de telles faces $F_1$ et $F_2$, il existe $i,j \in \{0,\cdots,n\}$,
$i \neq j$, tels que, d'une part, $u_i$ et $u_j$ sont des sommets de $F_1$,
d'autre part, $v_i$ et $v_j$ sont des sommets de $F_2$.
Les points $u_j+v_i$, $u_i+v_i$ et $u_j+v_j$ appartiennent alors \`a $F$ et v\'erifient
l'\'egalit\'e $u_j+v_i=2/3(u_i+v_i)+1/3(u_j+v_j)$.
L'application $\nu$ \'etant affine sur $F$, on en d\'eduit
que $\nu(u_j+v_i)=2/3\nu(u_i+v_i)+1/3\nu(u_j+v_j)$ i.e. $a_j+b_i=2/3(a_i+b_i)+1/3(a_j+b_j)$ 
et donc $2a_i-b_i=2a_j-b_j$.
{\cqfd}
\smallskip

Pour se fixer les id\'ees, on va maintenant supposer que la subdivision $\tau$ est la subdivision mixte
dont les $n$-polytopes sont les polytopes
$$F_l=F_1^l+F_2^{n-l} , \quad l=0,\cdots,n$$
o\`u l'on pose
$$
F_1^l=conv\{u_0,\cdots,u_l\}, \quad F_2^{n-l}=conv\{v_l,\cdots,v_n\}.$$

Pour cela il suffit de choisir des valeurs pour les $a_i$ et $b_i$
telles que $2a_0-b_0<2a_1-b_1<\cdots < 2a_n-b_n$.

\begin{ex}
\label{ex1}
Si $b_i=0$ et $ a_i=ki$
pour $i=0,\cdots,n$ de telle sorte que $\nu_2$ soit l'application nulle et $\nu_1$ la restriction de l'application
$(y_1,\cdots,y_n) \mapsto y_1+2y_2+\cdots+ny_n$, alors la subdivision $\tau$ associ\'ee est la subdivision mixte
dont les $n$-polytopes sont les polytopes
$$F_l=F_1^l+F_2^{n-l} , \quad l=0,\cdots,n.$$ 
\end{ex}

On note comme dans la section pr\'ec\'edente $f_i$ ($i=1,2$) le polyn\^ome obtenu en posant $t=1$ dans $f_{i,t}$,
puis on pose $f_1^l=f_1^{F_1^l}$
et $f_2^{n-l}=f_2^{F_2^{n-l}}$
(notons que $f_k^n=f_1$ et que $f_{2k}^n=f_2$).
On fait l'hypoth\`ese suivante.
\medskip

{\bf Hypoth\`ese.} $f_1$ et $f_2$ sont des polyn\^omes non d\'eg\'en\'er\'es.
\medskip

Par d\'efinition, chacun des $f_1^l$ et $f_2^{n-l}$ est alors non d\'eg\'en\'er\'e.
Les deux hypoth\`eses de la sous-section \ref{doubcons} sont donc maintenant satisfaites.
Rappelons que $X(F_1^l)=\{Z_{l+1}=\cdots=Z_n=0\} \subset {\C}P^n$
et que $X(F_2^{n-l})=\{Z_0=\cdots=Z_{l-1}=0\} \subset {\C}P^n$. Les polyn\^omes $f_1^l$
et $f_2^{n-l}$ d\'efinissent des hypersurfaces dans $X(F_1^l)$ et $X(F_2^{n-l})$,
que l'on note $X_k^l$ et $X_{2k}^{n-l}$, respectivement.
On a alors 
$$X_k^l=X_k^n \cap X(F_1^l) \quad \mbox{et} \quad X_{2k}^{n-l}=X_{2k}^n \cap X(F_2^{n-l}).$$
Chacun des $f_i^l$ \'etant non d\'eg\'en\'er\'es, chacune des hypersurfaces
${\R}X_k^l$ et ${\R}X_{2k}^{n-l}$ est non singuli\`ere.
\smallskip

Soit $l$ un nombre entier compris entre $1$ et $n-1$.
On consid\`ere la tranformation affine $\Delta_1$ (resp. $\Delta_2$, $\Delta$) permutant
$u_0$ et $u_l$ (resp. $v_0$ et $v_l$, $u_0+v_0$ et $u_l+v_l$) et laissant fixe chacun des autres sommets
de $T_1$ (resp. $T_2$, $T$). Les transformations $\Delta_1$, $\Delta_2$ et $\Delta$ appartiennent \`a $AFF_n(\Z)$
et ont pour m\^eme partie lin\'eaire la matrice $L \in GL_n(\Z)$ d\'efinie par $L_{ij}=-1$ si $i=l$ et $L_{ij}=\delta_{ij}$ sinon.
Jusqu'alors le syst\`eme de coordonn\'ees utilis\'e pour le tore alg\'ebrique
complexe ${({\C}^*)}^n={\C}P^n \setminus \cup_{i=0}^n\{Z_i=0\}$
\'etait le syst\`eme de coordonn\'ees $z=(z_1,\cdots,z_n)$, $z_i=Z_i/Z_0$, de la carte affine
$\{Z_0 \neq 0\}$. Les tranformations affines $\Delta_1$, $\Delta_2$ et $\Delta$ correspondent au m\^eme changement
de coordonn\'ees multiplicatif $z \mapsto \tilde{z}$ d\'efini par
$\tilde{z}=(\tilde{z_1},\cdots,\tilde{z_n})$ avec $\tilde{z}_l=Z_0/Z_l$ et $\tilde{z}_i=Z_i/Z_l$ si $i \neq l$.
Soient $\tilde{f}_1^l$ et $\tilde{f}_2^{n-l}$ les images dans ${\R}[\tilde{z}]$ de
$f_1^l$ et $f_2^{n-l}$ par $\Delta_1$ et $\Delta_2$, respectivement.
Le polytope de Newton de $\tilde{f}_1^l$ est le polytope $\Delta_1(F_1^l)$ et celui de
$\tilde{f}_2^{n-l}$ est $\Delta_2(F_2^{n-l})$.
On d\'ecompose l'espace ${\R}^n$ contenant ces polytopes comme le produit ${\R}^l \times {\R}^{n-l}$ en identifiant
${\R}^l$ \`a $\{(y_1,\cdots,y_n) \in {\R}^n \, | \, y_{l+1}= \cdots =y_n=0\}$ et
${\R}^{n-l}$ \`a $\{(y_1,\cdots,y_n) \in {\R}^n \, | \, y_{1}= \cdots =y_l=0\}$.
Le polyn\^ome $\tilde{f}_1^l$ est alors vu comme un polyn\^ome appartenant \`a $\R[\tilde{z_1},
\cdots,\tilde{z_l}]$ et ayant pour polytope de Newton $\Delta_1(F_1^l) \subset {\R}^l$.
De m\^eme, le polyn\^ome $\tilde{f}_2^{n-l}$ est alors vu comme un polyn\^ome appartenant \`a 
$\R[\tilde{z_{l+1}},\cdots,\tilde{z_n}]$ et ayant pour polytope de Newton $\Delta_2(F_2^{n-l}) \subset {\R}^{n-l}$.
Les hypersurfaces $X_k^l \subset X(F_1^l)$ et $X_{2k}^{n-l} \subset X(F_2^{n-l})$ sont d\'efinies
par les homog\'einis\'es de $\tilde{f}_1^l$ et de $\tilde{f}_2^{n-l}$ (qui coincident avec les homog\'einis\'es
de $f_1^l$ et de $f_2^{n-l}$), respectivement.

Notons $\pi_1$ (resp. $\pi_2$) la projection de ${({\Z}/2)}^n={({\Z}/2)}^l \times {({\Z}/2)}^{n-l}$
sur son facteur ${({\Z}/2)}^l$ (resp. ${({\Z}/2)}^{n-l}$) et $\pi:{({\Z}/2)}^n \rightarrow {({\Z}/2)}^l \times {({\Z}/2)}^{n-l}$
l'application $g \mapsto (\pi_1(g),\pi_2(g))$.
Pour tout $g \in {({\Z}/2)}^n$, on a 
alors $g(\Delta_1(F_1^l))=g_1(\Delta_1(F_1^l))$, $g(\Delta_2(F_2^{n-l}))=g_2(\Delta_2(F_2^{n-l}))$ si $\pi(g)=(g_1,g_2)$.

D\'esignons par $c.c \; {({\R}^*)}^n$ la r\'eunion des $2^n$ composantes connexes (que l'on appellera orthants)
de ${({\R}^*)}^n={\R}P^n \setminus \cup_{i=0}^n \{Z_i=0\}$.
Soient $\psi_0$ et $\psi_l$ les param\'etrisations (bijections) ${({\Z}/2)}^n \rightarrow c.c \; {({\R}^*)}^n$
d\'efinies par
$\psi_0(g)=\{(z_1,\cdots,z_n) \in {\R}^n \; | \; {(-1)}^{\epsilon_i}z_i>0, \; i=1,\cdots,n\}$
et $\psi_l(g)=\{(\tilde{z}_1,\cdots,\tilde{z}_n) \in {\R}^n \; | \; {(-1)}^{\epsilon_i}\tilde{z}_i>0, \; i=1,\cdots,n\}$,
si $g=(\epsilon_1,\cdots,\epsilon_n)$ . On a $\psi_0(g)=\R(g)$ dans les notations pr\'ec\'edentes.
On a d\'eja d\'efini l'orthant positif comme \'etant ${({\R}_+^*)}^n=\{(z_1,\cdots,z_n) \in {\R}^n \; | \; z_i>0, \; i=1,\cdots,n\}$.
Notons que l'on a ${({\R}_+^*)}^n=\psi_0(id)=\psi_l(id)$, puis que  $\psi_0(g)=g({({\R}_+^*)}^n)$ (resp. $\psi_l(g)=g({({\R}_+^*)}^n)$)
pour tout $g \in{({\Z}/2)}^n$ si l'on identifie \`a ${({\Z}/2)}^n$ le groupe des sym\'etries de
${\R}^n$ par rapport aux hyperplans de coordonn\'ees
$\{{z}_i=0\}$ (resp, $\{\tilde{z}_i=0\}$) comme dans la sous-section \ref{carte}.

Soit $c.c {({\R}^*)}^l$ (resp. $c.c {({\R}^*)}^{n-l}$) l'union des
composantes connexes de ${({\R}^*)}^l=\R X(F_1^l) \setminus \cup_{i=0}^l \{Z_i=0\}$
(resp. ${({\R}^*)}^{n-l}=\R X(F_2^{n-l}) \setminus \cup_{i=l}^n \{Z_i=0\}$). Par restriction,
on obtient des param\'etrisations $\psi_{l,1}:{({\Z}/2)}^l \rightarrow c.c {({\R}^*)}^l$,
et $\psi_{l,2}:{({\Z}/2)}^{n-l} \rightarrow c.c {({\R}^*)}^{n-l}$ telles que
$\psi_l(g)=\psi_{l,1}(g_1) \times \psi_{l,2}(g_2)$ pour tout $g \in {({\Z}/2)}^n$ avec $\pi(g)=(g_1,g_2)$.
Soit $\alpha_l=\psi_l^{-1} \circ \psi_0:{({\Z}/2)}^n \rightarrow {({\Z}/2)}^n$ l'application de changement de param\'etrisation.
On a alors $\R(g)=\psi_0(g)=\psi_l(\alpha_l(g))$ pour tout $g \in {({\Z}/2)}^n$.
\smallskip

Pour tout nombre entier $l$ compris entre $0$ et $n$, on note
$$Y_{2k}^{n-l}=\{U^2-{\bar f}_2^{n-l}(Z)=0\} \subset {\C}P^{n-l+1}(1,k)$$
le plan doubl\'e associ\'e \`a ${\bar f}_2^{n-l}$
et
$${\R}P^{n-l}_+=\{{\bar f}_2^{n-l} \geq 0\} \subset {\R}X(F_2^{n-l}).$$
Rappelons que ${\R}Y_{2k}^{n-l}$ est non singulier ($f_2^{n-l}$ est non d\'eg\'en\'er\'e),
se projette sur ${\R}P^{n-l}_+$ et se ramifie sur ${\R}X_{2k}^{n-l}$.
On note, respectivement,
$${\R}X_k^l(\star) \; , \quad {\R}X_{2k}^{n-l}(\star) \; , \quad {\R}P^{n-l}_+(\star) \quad \mbox{et}  \quad {\R}Y_{2k}^{n-l}(\star)$$
les intersections ${\R}X_k^l \cap {({\R}^*)}^l$, ${\R}X_{2k}^{n-l} \cap {({\R}^*)}^{n-l}$,
${\R}P^{n-l}_+ \cap {({\R}^*)}^{n-l}$ et la partie de ${\R}Y_{2k}^{n-l}$  se projettant sur ${\R}P^{n-l}_+(\star)$.
%

De la m\^eme mani\`ere, si $1 \leq l \leq n-1$ et si $(g_1,g_2) \in {({\Z}/2)}^l \times {({\Z}/2)}^{n-l}$, on note, respectivement,
$${\R}X_k^l(g_1) \; , \quad {\R}X_{2k}^{n-l}(g_2) \; , \quad {\R}P^{n-l}_+(g_2) \quad \mbox{et}  \quad {\R}Y_{2k}^{n-l}(g_2)$$
les intersections ${\R}X_k^l \cap \psi_{l,1}(g_1)$, ${\R}X_{2k}^{n-l} \cap \psi_{l,2}(g_2)$,
${\R}P^{n-l}_+ \cap \psi_{l,2}(g_2)$ et la partie de ${\R}Y_{2k}^{n-l}$  se projettant sur ${\R}P^{n-l}_+(g_2)$.

\begin{lem}
\label{cha}
Soit $g \in {({\Z}/2)}^n$.

Pour tout nombre entier $l$ compris entre $1$ et $n-1$,
si $(g_1,g_2)=(\pi \circ \alpha_l)(g)$, alors on a les hom\'eomorphismes suivants:
$$C_{F_1^l}^g(f_1^l) \simeq {\R}X_k^l(g_1) \; , \quad
C_{F_2^{n-l}}^g(f_2^{n-l})  \simeq {\R}X_{2k}(g_2)$$
$$
C_{F_2^{n-l},+}^g(f_2^{n-l}) \simeq {\R}P^{n-l}_+(g_2)$$
\end{lem}
{\sc Preuve.}
Soient $\phi:{({\C}^*)}^n \rightarrow F_1^l$ l'application moment associ\'ee \`a $F_1^l$ et
$\phi_l:{({\C}^*)}^l \rightarrow \Delta_1(F_1^l)$ l'application moment associ\'ee \`a $\Delta_1(F_1^l)$,
relativement aux syst\`emes de cooordonn\'ees $(z_1,\cdots,z_n)$ pour ${({\C}^*)}^n$
et $(\tilde{z}_1,\cdots,\tilde{z}_l)$ pour ${({\C}^*)}^l$. Sachant que $z^w=\tilde{z}_l^{-k}{\tilde{z}}^{\Delta_1(w)}$,
on en d\'eduit que pour tout $(z_1,\cdots,z_n)$ on a 
$(\Delta_1 \circ \phi)(z_1,\cdots,z_n)=\phi_l(\tilde{z}_1,\cdots,\tilde{z}_l)$,
puis que
$$C_{F_1^l}^{id}(f_1^l) \stackrel{\Delta_1}{\simeq} C_{\Delta_1(F_1^l)}^{id}(\tilde{f}_1^l).$$
Soit $g \in {({\Z}/2)}^n$ et $g_1=(\pi_1 \circ \alpha_l)(g)$.
Sachant que $z=(z_1,\cdots,z_n) \in {\R}(g)=\psi_0(g) \Rightarrow (\tilde{z}_1,\cdots,\tilde{z}_l) \in \psi_{l,1}(g_1)$, on obtient que
$$C_{F_1^l}^{id}(f_1^l \circ g) \stackrel{\Delta_1}{\simeq} C_{\Delta_1(F_1^l)}^{id}(\tilde{f}_1^l \circ g_1).$$
Maintenant, par d\'efinition on a
$$C_{F_1^l}^{g}(f_1^l)=g \left(C_{F_1^l}^{id}(f_1^l \circ g)\right) \; , \quad
C_{\Delta_1(F_1^l)}^{g_1}(\tilde{f}_1^l) =g_1 \left(C_{\Delta_1(F_1^l)}^{id}(\tilde{f}_1^l \circ g_1) \right).$$
Par cons\'equent $g_1 \circ \Delta_1 \circ g:g(F_1^l) \rightarrow (g_1 \circ \Delta_1)(F_1^l))$
est un hom\'eomorphisme envoyant $C_{F_1^l}^{g}(f_1^l)$ sur $C_{\Delta_1(F_1^l)}^{g_1}(\tilde{f}_1^l)$.
Pour finir, le polytope $\Delta_1(F_1^l)$ \'etant de dimension $l$, la restriction
de $\phi_l$ \`a ${({\R}_+^*)}^l$ est un hom\'eomorphisme sur l'int\'erieur de $\Delta_1(F_1^l)$,
et donc $g_1 \circ  \phi_l \circ g_1:\psi_{l,1}(g_1) \rightarrow (g_1 \circ \Delta_1)(F_1^l)$ est un hom\'eomorphisme
envoyant ${\R}X_k^l(g_1)$ sur $C_{\Delta_1(F_1^l)}^{g_1}(\tilde{f}_1^l)$.
La d\'emonstration de l'existence des deux autres hom\'eomorphismes est tout \`a fait similaire.
{\cqfd}

La remarque qui suit nous sera utile dans la d\'emonstration de la proposition plus bas.

\begin{rem}
\label{extr}

Sachant que ${\R}X_k^0={\R}X_{2k}^0=\emptyset$, pour tout $g \in {({\Z}/2)}^n$ on a
$$
C_{F_1^0}^g(f_1^0)=C_{F_2^0}^g(f_{2}^0)=\emptyset.$$

Rappelons que $F_2^0=(0,\cdots,2k)$ et que $f_{2}^0=f_{2}^{F_2^0}=a_{2,(0,\cdots,0,2k)} z_n^{2k}$.
On obtient que pour tout $g \in {({\Z}/2)}^n$ on a

\begin{itemize}
\item $a_{2,(0,\cdots,0,2k)}>0 \Rightarrow$
$$C_{F_2^0,+}^g(f_{2}^0)=C_{F_2^0}=\{(0,\cdots,0,2k)\},
\quad {\R}P^{0}_+=X(F_2^0)=\{(0:\cdots:0:1)\},$$
\item $a_{2,(0,\cdots,0,2k)}<0 \Rightarrow$ 
$$C_{F_2^0,+}^g(f_{2}^0)=\emptyset= {\R}P^{0}_+.$$
\end{itemize}
\end{rem}

On peut maintenant \'enoncer le r\'esultat principal de cette section.

\begin{prop}
\label{doub}
Pour $t,\epsilon>0$ suffisamment petits l'hypersurface doubl\'ee ${\R}X$ est non singuli\`ere
et hom\'eomorphe
au collage d\'etermin\'e par la subdivision $\tau$ des vari\'et\'es produit
$${\R}X_k^l(\star) \times {\R}Y_{2k}^{n-l}(\star), \quad l=1,\cdots,n.$$
\end{prop}
{\sc Preuve.}
Soit $Y=\{U^2-{\bar f}_{2,t}(Z)=0,{\bar f}_{1,t}(Z)=0\} \subset {\C}P^{n+1}(1,k)$, $t>0$ petit,
l'hypersurface r\'eelle se d\'eformant sur $X$ (${\R}Y$ est hom\'eomorphe \`a ${\R}X$).
On note comme dans la sous-section \ref{doubcons}
$X_k$ et $X_{2k}$ les hypersurfaces d\'efinies pour $t>0$ petit par 
$f_{1,t}$ et $f_{2,t}$, respectivement, et on pose ${\R}X_{k,+}:=
\{Z \in {\R}X_k, \; {\bar f}_{2,t}(Z) \geq 0\}$, ${\R}X_{k,+}(g):={\R}X_{k,+}
\cap {\R}(g)$.

Soit $g \in {({\Z}/2)}^n$. On montre que
la partie ${\R}Y(g)$ de ${\R}Y$ se projettant sur ${\R}X_{k,+}(g)$
est hom\'eomorphe au collage des
$${\R}X_k^l(g_1) \times {\R}Y_{2k}^{n-l}(g_2), \quad l=1,\cdots,n-1,$$
tels que $(g_1,g_2)=(\pi \circ \alpha_l)(g)$ et de
$$\big({\R}X_k^n \cap \R(g)\big) \times {\R}Y_{2k}^0.$$
D'apr\`es le corollaire \ref{hom}, l'espace ${\R}X_{k,+}(g)$
est hom\'eomorphe \`a $S_2^+ \cap S_1 \cap {\R}(g)$ qui s'obtient par collage des
$$C_{F_1^l}^{g}(f_1^l)+C_{F_2^{n-l},+}^{g}(f_2^{n-l}) \subset
g(F_1^l)+g(F_2^{n-l})=g(F^l), \quad l=0,\cdots,n.$$
En utilisant Le Lemme \ref{cha}, la remarque \ref{extr} et le fait que $C_{F_1^n}^{g}(f_1^n)$
est hom\'eomorphe \`a ${\R}X_k^n \cap \R(g)$, on obtient alors que ${\R}X_{k,+}(g)$ 
est hom\'eomorphe au collage des
$${\R}X_k^l(g_1) \times {\R}P^{n-l}_+(g_2), \quad l=1,\cdots,n-1 $$
tels que $(g_1, g_2)=(\pi \circ \alpha_l)(g)$ et de
$$\big({\R}X_k^n \cap \R(g)\big) \times {\R}P_+^0.$$

Sachant que ${\R}Y(g)$ est hom\'eomorphe au double de ${\R}X_{k,+}(g)$, on obtient que
${\R}Y(g)$ est hom\'eomorphe au collage pour $l=1,\cdots,n-1$ des doubles
des ${\R}X_k^l(g_1) \times {\R}P^{n-l}_+(g_2)$
v\'erifiant $(g_1, g_2)=(\pi \circ \alpha_l)(g)$
et du double de $\big({\R}X_k^n \cap \R(g)\big) \times {\R}P_+^0$.
Il reste \`a remarquer que le double de ${\R}X_k^l(g_1) \times {\R}P^{n-l}_+(g_2)$
est hom\'eomorphe \`a ${\R}X_k^l(g_1) \times {\R}Y_{2k}^{n-l}(g_2)$
et que le double de $\big({\R}X_k^n \cap \R(g)\big) \times {\R}P_+^0$
est hom\'eomorphe \`a $\big({\R}X_k^n \cap \R(g)\big) \times {\R}Y_{2k}^0$.
{\cqfd}

Notons que ${\R}Y_{2k}^0$ est soit vide, soit
r\'eduit \`a deux points distincts suivant le signe du coefficient devant $z_n^{2k}$
dans $f_2(z)$ (voir la remarque \ref{extr}).
\smallskip

Rappelons que $X_k$ (resp. $X_{2k}$) d\'esigne l'hypersurface de ${\C}P^n$
d\'efinie par $f_{1,t}$ (resp. $f_{2.t}$) pour $t>0$ suffisamment petit.
La remarque qui suit, qui ne sera pas utilis\'ee, implique que dans la proposition \ref{doub}
on peut remplacer $X_k^l$ par $X_k \cap X(F_1^l)$ et $Y_{2k}^{n-l}$ par le plan doubl\'e associ\'e au
tronqu\'e de $f_{2,t}$ sur $F_2^{n-l}$.

\begin{rem}
$\nu_i$ \'etant la restriction d'une application affine
sur ${\R}^n$, on a, pour tout $0<t \leq 1$ et tout $l=0,\cdots,n$, des hom\'eomorphismes de paires
$$({\R}X(F_1^l), {\R}X_k^l) \simeq ({\R}X(F_1^l), {\R}X(F_1^l) \cap {\R}X_k),$$
$$({\R}X(F_2^{n-l}), {\R}X_{2k}^{n-l}) \simeq ({\R}X(F_2^{n-l}),{\R}X(F_2^{n-l})
\cap {\R}X_{2k}).$$
\end{rem}
Cette remarque est une cons\'equence directe de l'observation suivante. 
Soit $f_t(z)=\sum_{w \in W} a_wt^{\nu(w)}z^w$ un polyn\^ome de Viro associ\'e \`a une application
$\nu$ restriction d'une application affine
$(y_1,\cdots,y_n) \mapsto \sum_{i=1}^n{\lambda_iy_i}+\lambda_0$.
Si $f$ d\'esigne le polyn\^ome obtenu en posant $t=1$ dans $f_t$,
alors $f_t(z)=t^{\lambda_0}f(\tilde{z})$ o\`u
$z=(z_1,\cdots,z_n) \mapsto \tilde{z}=(\tilde{z}_1,\cdots,\tilde{z}_n)$
est le changement de coordonn\'ees du tore ${({\C}^*)}^n$ d\'efini par $\tilde{z_i}=t^{\lambda_i}z_i$.

\begin{ex}[Suite de l'exemple \ref{ex1}]
Si $\nu_1$ est la restriction de l'application $(y_1,\cdots,y_n) \mapsto y_1+2y_2+\cdots+ny_n$
et $\nu_2$ est l'application nulle, alors
les hom\'eomorphismes ${\R}X(F_1^l) \rightarrow {\R}X(F_1^l)$ plus haut
sont donn\'es par les restrictions
de l'application $(Z_0:Z_1:Z_2:\cdots:Z_n) \mapsto (Z_0:tZ_1:t^2Z_2:\cdots:t^nZ_n)$ et les hom\'eomorphismes
${\R}X(F_2^{n-l}) \rightarrow {\R}X(F_2^{n-l})$ sont donn\'es par l'identit\'e.
\end{ex}
%
%
%

\section{Asymptotiques de nombres de Betti}

\subsection{Comportement asymptotique}

Soit $i$ un entier positif.

On note, respectivement,
$$Max \; b_i({\R}  X_m^n) \quad \mbox{et} \quad Max \; b_i({\R}Y_{2k}^n) $$
la valeur maximale
des nombres de Betti $b_i({\R}  X_m^n)$ prise sur l'ensemble des hypersurfaces $X_m^n$
de degr\'e $m$ dans ${\C}P^n$ pour $m$ et $n$ fix\'es, et la valeur maximale des nombres de Betti $b_i({\R}Y_{2k}^n)$
prise sur l'ensemble des plans doubl\'es $Y_{2k}^n=\{U^2-{\bar{f}}_{2k}(Z)=0\} \subset {\C}P^{n+1}(1,k)$ associ\'es
\`a des polyn\^omes $f_{2k}$ de degr\'e $2k$ pour $k$ et $n$ fix\'es.

On s'int\'eresse aux suites
$$\big(Max \; b_i({\R}X_m^n)\big)_{m \geq 1} \; , \quad \big(Max \; b_i({\R}Y_{2k}^n)\big)_{k\geq 1},$$
et plus particuli\`erement \`a leurs comportements asymptotiques.

Ici et par la suite, l'expression
$${\mathcal R}(n,m)$$
d\'esignera une fonction ``reste'' de $n$ et de $m$ comprise entre deux fonctions polynomiales de degr\'e $n$ en $m$
(dont les coefficients ne d\'ependent pas de $m$). En particulier, on aura
$${\mathcal R}(n,m)/m^{n+1} \rightarrow 0 \quad \mbox{pour} \quad m \rightarrow +\infty.$$

L'in\'egalit\'e de Smith-Thom (voir l'introduction) implique imm\'ediatement
que $Max \; b_i({\R}  X_m^n) \leq b_*(X_m^n)$ et $ Max \; b_i({\R}Y_{2k}^n) \leq  b_*(Y_{2k}^n)$.
\`A titre d'information, en utilisant \cite{Hodge} (on pourra voir aussi  \cite{Kha} section 2.2), on obtient les formules suivantes
$$
\begin{array}{l}
b_*(X_m^n)=\frac{(m-1)^{n+1}+{(-1)}^n}{m}+n+{(-1)}^{n+1} \\
  \\
b_*(Y_{2k}^n)=\frac{(2k-1)^{n+1}+{(-1)}^n}{2k}+n+1.
\end{array}
$$
\smallskip

Cela donne en particulier

$$Max \; b_i({\R}X_m^n) \leq m^n + {\mathcal R}(n-1,m), \quad Max \; b_i({\R}Y_{2k}^n) \leq 2^n k^n+{\mathcal R}(n-1,k).$$

\begin{prop}
\label{as1}
Il existe des nombres r\'eels ${\h}_{i,n}$ et ${\d}_{i,n}$
pour lesquels on ait les \'equivalences asymptotiques suivantes
$$Max \; b_i({\R}X_m^n) \sim {\h}_{i,n} \cdot m^n, \quad Max \;
b_i({\R}Y_{2k}^n) \sim {\d}_{i,n} \cdot k^n$$
pour $m$ (resp. $k$) $\rightarrow + \infty$.
\end{prop}

{\sc preuve.}
On commence par d\'emontrer l'existence de ${\h}_{i,n}$. La d\'emonstration qui suit est
une g\'en\'eralisation directe de celle donn\'ee dans \cite{Basymp} de l'existence de ${\h}_{i,3}$.
L'in\'egalit\'e de Smith-Thom implique que la suite ($Max \; b_i({\R}X_m^n)/m^n)_{m \geq 1}$
est born\'ee, et donc admet une limite sup\'erieure $L$.
On montre que cette suite converge vers $L$ de la mani\`ere suivante: pour un $\epsilon>0$ donn\'e, on construit
avec la m\'ethode de viro une famille d'hypersurfaces non singuli\`eres ${\R}X_m^n$
telle que $b_i({\R}X_m^n)/m^n>L-\epsilon$ pour tout degr\'e $m$
suffisamment grand.

Soit $\epsilon>0$. Par d\'efinition de $L$, il existe
un hypersurface $X_{m_0}^n$ de suffisamment grand degr\'e ${m_0}$ telle que
$b_i({\R}X_{m_0}^n)/{m_0}^n>L-\epsilon/2$. On peut choisir $X_{m_0}^n$ de telle sorte
que ${\R}X_{m_0}^n$ n'intersecte aucun des hyperplans de coordonn\'ees $\{Z_i=0\}$, autrement dit que
${\R}X_{m_0}^n$ est contenu dans le tore r\'eel ${({\R}^*)}^n$ (en particulier $m_0$ doit \^etre pair).

Rappelons que dans ce papier par polytope (ou simplexe ...) on entend polytope
convexe \`a sommets entiers. De m\^eme, par subdivision ou triangulation d'un polytope, on entendra
subdivision ou triangulation \`a sommets entiers. Rappelons \'egalement que l'on note
$T_m^n$ le polytope de Newton d'un polyn\^ome affine
g\'en\'erique de degr\'e $m$ en $n$ variables.

Pour tout nombre entier $p \geq 1$, on choisit une triangulation convexe de $T_p^n$
dont tous les $n$-simplexes sont
de volume euclidien $1/n!$ ou autrement dit qui contient $p^n$ $n$-simplexes
(une telle triangulation est dite primitive, ou unimodulaire, on peut voir \cite{ItVi} pour des exemples).
Chacun de ces $n$-simplexes est l'image de $T_1^n$ par une transformation de $AFF_n(\Z)$. 
Consid\'erons maintenant la triangulation de $T_{p({m_0}+n+1)}^n$
obtenue en appliquant l'homoth\'etie de rapport ${m_0}+n+1$ \`a la pr\'ec\'edente triangulation.
La triangulation de $T_{p({m_0}+n+1)}^n$ que l'on obtient est bien sur convexe et est constitu\'ee
de $p^n$ $n$-simplexes de volume euclidien $({m_0}+n+1)^n/n!$. Chacun de ces $n$-simplexes est
l'image de $T_{{m_0}+n+1}^n$ par une transformation de $AFF_n(\Z)$.
Sachant que $(1,\cdots,1)+T_{m_0}^n$ est contenu \`a l'int\'erieur de $T_{{m_0}+n+1}^n$,
on peut donc placer \`a l'int\'erieur de chacun des $n$-simplexes de la triangulation de $T_{p({m_0}+n+1)}^n$
l'image de $T_{m_0}^n$ par un \'el\'ement de $AFF_n(\Z)$.
On raffine alors la triangulation de $T_{p({m_0}+n+1)}^n$ en une triangulation convexe, que l'on note $\tau$,
contenant $p^n$ images disjointes de $T_{{m_0}}^n$ par des \'el\'ements de $AFF_n(\Z)$.

Les $p^n$ images de $T_{{m_0}}^n$ contenues dans $\tau$ \'etant disjointes, on peut trouver
un polyn\^ome $f$ de polytope de Newton $T_{p({m_0}+n+1)}^n$ v\'erifiant
les conditions suivantes.
\begin{itemize}
\item Si $F$ est l'une des $p^n$ images de $T_{{m_0}}^n$ mentionn\'ee plus haut,
alors $f^F$ est l'image d'un polyn\^ome non d\'eg\'en\'er\'e d\'efinissant l'hypersurface $X_{m_0}^n$
par la transformation de $AFF_n(\Z)$ correspondante.
\item $f^F$ est non d\'eg\'en\'er\'e pour tout $F \in \tau$.
\end{itemize}

Soit $f_t$ un polyn\^ome de Viro tel que $f_t=f$ pour $t=1$
et d\'efini par une fonction $\nu$
d\'eterminant la subdivision $\tau$. Le th\'eor\`eme de Viro implique alors que pour $t>0$
suffisamment petit l'hypersurface ${\R}X_{p({m_0}+n+1)}^n$
est non singuli\`ere et hom\'eomorphe au collage, d\'etermin\'e par $\tau$,
des hypersurfaces $\{f^F=0\} \cap {({\R}^*)}^n$, $F \in \tau$.
Si $F \in \tau$ est l'une des $p^n$ images de $T_{{m_0}}^n$ pr\'ec\'edentes, alors
$\{f^F=0\} \cap {({\R}^*)}^n $ est hom\'eomorphe \`a ${\R}X_{m_0}^n \cap {({\R}^*)}^n$,
et donc \`a ${\R}X_{m_0}^n$ puisque, par hypoth\`ese, ${\R}X_{m_0}^n$ n'intersecte
pas les hyperplans de coordonn\'ees de ${\C}P^n$. On obtient donc $p^n$ copies hom\'eomorphes de ${\R}X_{m_0}^n$
qui sont disjointes dans ${\R}X_{p({m_0}+n+1)}^n$.
Par cons\'equent, $b_i({\R}X_{p({m_0}+n+1)}) \geq p^n b_i({\R}X_{m_0}^n)$, et donc
\begin{equation}
\label{equ1}
\frac{b_i({\R}X_{p({m_0}+n+1)})}{p^n \cdot ({m_0}+n+1)^n}
 \geq \frac{b_i({\R}X_{m_0}^n)}{{m_0}^n} - \frac{b_i({\R}X_{m_0}^n)}{{m_0}^n}\cdot \big[1-\frac{1}{(1+(n+1)/{m_0})^n}\big].
\end{equation}
Sachant que $\frac{b_i({\R}X_{m_0}^n)}{{m_0}^n} \geq L-\epsilon/2$, et que
la suite $(Max \; b_i({\R}X_m^n)/m^n)_{m \geq 1}$ est born\'ee,
l'in\'egalit\'e (\ref{equ1}) implique que, quitte \`a choisir une surface $X_{m_0}^n$
de plus grand degr\'e ${m_0}$, on a 
\begin{equation}
\label{equ2}
\frac{b_i({\R}X_{p({m_0}+n+1)})}{p^n \cdot ({m_0}+n+1)^n}
\geq \frac{b_i({\R}X_{m_0}^n)}{{m_0}^n}-\epsilon/4.
\end{equation}

Soit $p_0 \geq 1$ un nombre entier et soit $\{X_{p({m_0}+n+1)^n}, \; p \geq p_0\}$ une famille form\'ee
d'hypersurfaces construites comme pr\'ec\'edemment pour chaque $p \geq p_0$ et v\'erifiant (\ref{equ2}).
On compl\`ete cette famille en une famille ${\mathcal F}=\{X_m^n, \; m \geq p_0({m_0}+n+1) \}$ de la mani\`ere suivante
Pour tout nombre entier $m$ tel que $p({m_0}+n+1)< m < (p+1)({m_0}+n+1)$ avec $p \geq p_0$,
on construit une hypersurface $X_m^n$ en lissifiant
l'union de $X_{p({m_0}+n+1)}$ avec $m-p({m_0}+n+1)$ hyperplans n'intersectant pas les
$p^n$ copies hom\'eomorphes de ${\R}X_{m_0}^n$ qui sont contenues dans ${\R}X_{p({m_0}+n+1)}^n$
(on peut prendre des hyperplans proches des hyperplans de coordonn\'ees).
On a alors
$b_i({\R}X_m^n)\geq p^n \cdot b_i({\R}X_{m_0}^n)$, et donc, sachant que $m < (p+1)({m_0}+n+1)$
\begin{equation}
\label{equ3}
\frac{b_i({\R}X_m^n)}{m^n}\geq \frac{b_i({\R}X_{m_0}^n)}{{m_0}^n} - \frac{b_i({\R}X_{m_0}^n)}{{m_0}^n}\cdot \big[1-\frac{1}{[(1+(n+1)/{m_0})(1+1/p)]^n}\big].
\end{equation}
Maintenant, les in\'egalit\'es (\ref{equ2}), (\ref{equ3})
et le fait que $(Max \; b_i({\R}X_m^n)/m^n)_{m \geq 1}$ soit born\'ee
impliquent que si $p_0$ est suffisamment grand, alors toute hypersurface appartenant \`a  ${\mathcal F}$
v\'erifie
\begin{equation}
\label{equ4}
\frac{b_i({\R}X_m^n)}{m^n}
\geq \frac{b_i({\R}X_{m_0}^n)}{{m_0}^n}-\epsilon/2.
\end{equation}
ce qui, avec $\frac{b_i({\R}X_{m_0}^n)}{{m_0}^n} \geq L-\epsilon/2$, implique que
$$
\frac{b_i({\R}X_m^n)}{m^n} \geq L -\epsilon.
$$ 
Cela termine la d\'emonstration de l'existence de $\h_{i,n}$.

L'existence de $\d_{i,n}$ peut se montrer exactement de la m\^eme mani\`ere.
Soit $L$ la limite sup\'erieure de la suite $Max \; b_i({\R}Y_{2k}^n)/k^n$,
et soit $\epsilon>0$. Il existe
un polyn\^ome r\'eel $f_{2k_0}$ non d\'eg\'en\'er\'e (de degr\'e $2k_0$ en $n$ variables)
tel que le plan doubl\'e r\'eel associ\'e $Y_{2k_0}^n=\{U^2-{\bar f}_{2k_0}(Z)\} \subset {\C}P^{n+1}(1,k_0)$
v\'erifie $b_i({\R}Y_{2k_0}^n)/{k_0}^n \geq L-\epsilon/2$.
On peut de plus, quitte \`a augmenter $k_0$, choisir ce polyn\^ome de telle sorte que $\{{\bar f}_{2k_0} \geq 0\}$
n'intersecte pas les hyperplans de coordonn\'ees de ${\R}P^n$.
Pour ce faire, on peut utiliser la m\'ethode de Viro et consid\'erer une
subdivision poly\'edrale convexe de $T_{2(k_0+n+1)}^n$ qui contienne
le simplexe $(2,\cdots,2)+T_{2k_0}^n$ (ce dernier est situ\'e \`a l'int\'erieur de $T_{2(k_0+n+1)}^n$),
ainsi qu'un polyn\^ome $f$ de polytope de Newton $T_{2(k_0+n+1)}^n$
v\'erifiant les conditions suivantes:
\begin{itemize}
\item le tronqu\'e de $f$ sur $(2,\cdots,2)+T_{2k_0}^n$ est $(z_1 \cdots z_n)^2 f_{2k_0}(z)$,
\item si $F$ est un polytope de la subdivision enti\`erement contenu dans une face de $T_{2(k_0+n+1)}^n$,
alors $\{f^F=0\} \cap {({\R}^*)}^n=\emptyset$,
\item pour tout polytope $F$ de la subdivision, le polyn\^ome $f^F$
est non d\'eg\'en\'er\'e.
\end{itemize}
Si $f_{2(k_0+n+1)}$ d\'esigne un polyn\^ome non d\'eg\'en\'er\'e obtenu par la m\'ethode de Viro
\`a partir de ces donn\'ees (obtenu pour $t>0$
suffisamment petit, \`a partir d'un polyn\^ome de Viro $f_t$ v\'erifiant $f_t=f$ pour $t=1$
et associ\'e \`a une fonction $\nu$ certifiant la convexit\'e de la subdivision pr\'ec\'dente),
alors clairement $\{f_{2(k_0+n+1)} \geq 0\}$
n'intersecte pas les hyperplans de coordonn\'ees de ${\R}P^n$, et de plus, si $k_0$
est suffisamment grand, le plan doubl\'e associ\'e $Y_{2(k_0+n+1)}^n=\{U^2-{\bar f}_{2(k_0+n+1)}(Z)\} \subset {\C}P^{n+1}(1,k_0)$
v\'erifie $b_i({\R}Y_{2(k_0+n+1)}^n)/{(k_0+n+1)}^n \geq L-\epsilon/2$.

Etant donn\'e un polyn\^ome r\'eel $f_{2k_0}$ non d\'eg\'en\'er\'e
tel que le plan doubl\'e r\'eel associ\'e $Y_{2k_0}^n=\{U^2-{\bar f}_{2k_0}(Z)\} \subset {\C}P^{n+1}(1,k_0)$
v\'erifie $b_i({\R}Y_{2k_0}^n)/{k_0}^n \geq L-\epsilon/2$ et tel que $\{{\bar f}_{2k_0} \geq 0\}$
n'intersecte pas les hyperplans de coordonn\'ees de ${\R}P^n$, on construit ensuite \`a l'aide de la m\'ethode
de Viro, pour tout entier $p \geq 1$, un polyn\^ome $f_{2p(k_0+n+1)}$ non d\'eg\'en\'er\'e de degr\'e $2p(k_0+n+1)$
tel que $\{{\bar f}_ {2p(k_0+n+1)}\geq 0\}$ contienne $p^n$ copies hom\'eomorphes et disjointes de $\{{\bar f}_{2k_0} \geq 0\}$.
On termine alors de la m\^eme mani\`ere que dans la d\'emonstration de l'existence de $\h_{i,n}$.
{\cqfd}

\subsection{Estimations des $\bf {\h}_{i,n}$ et $\bf {\d}_{i,n}$}
\label{est}
Les valeurs exactes des ${\h}_{i,n}$ et ${\d}_{i,n}$ sont connues
pour les premi\`eres dimensions (on a bien sur ${\h}_{i,n}=0$ si $i \geq n$
et ${\d}_{i,n}=0$ si $i \geq n+1$):
$$
{\d}_{0,0}=2,\quad
{\h}_{0,1}={\d}_{0,1}={\d}_{1,1}=1,\quad
{\h}_{0,2}={\h}_{1,2}=\frac{1}{2}.
$$

La valeur $\frac{1}{2}$ pour ${\h}_{0,2}$ et ${\h}_{1,2}$ est une cons\'equence du fameux th\'eor\`eme d'Harnack
disant que le nombre de composantes connexes de la partie r\'eelle d'une courbe
alg\'ebrique r\'eelle non singuli\`ere de degr\'e $m$ dans ${\C}P^2$ est major\'e
par $\frac{(m-1)(m-2)}{2}+1$ et que cette borne est exacte pour tout degr\'e.

On connait de plus les estimations suivantes \cite{Itrag}
$$\frac{27}{16}
\leq {\d}_{0,2} \leq \frac{7}{4}, \quad \frac{27}{8}
\leq {\d}_{1,2} \leq \frac{7}{2},$$
ainsi que \cite{Basymp}
$$\frac{13}{36}
\leq {\h}_{0,3} \leq \frac{5}{12}, \quad \frac{13}{18}
\leq {\h}_{1,3} \leq \frac{5}{6}.$$
Des bornes sup\'erieures pour les ${\h}_{i,n}$ et ${\d}_{i,n}$
sont obtenues de mani\`ere classique en utilisant les in\'egalit\'es de Smith-Thom
et de Comessatti-Petrowsky-Oleinik g\'en\'eralis\'ees.

\begin{prop}
\label{bsup}
Soient $c_n$ et $c_n'$ les coefficients dominants des nombres
de Hodge centraux des $X_m^n$ et $Y_{2k}^n$:
$$h^{\frac{n-1}{2},\frac{n-1}{2}}(X_m^n)=c_n m^n +{\mathcal R}(n-1,m),\quad
n \; \mbox{impair},$$
$$h^{\frac{n}{2},\frac{n}{2}}(Y_{2k}^n)=c_n' k^n +{\mathcal R}(n-1,k), \quad n \; \mbox{pair}.$$
\begin{enumerate}
\item Si $n$ est pair on a ${\h}_{i,n} \leq \frac{1}{2}$,
sinon on a ${\h}_{\frac{n-1}{2},n} \leq \frac{1+c_n}{2}$
et ${\h}_{i,n} \leq \frac{1+c_n}{4}$ pour $i \neq \frac{n-1}{2}$.

\item Si $n$ est impair on a ${\d}_{i,n} \leq 2^{n-1}$, sinon
on a ${\d}_{\frac{n}{2},n} \leq \frac{2^n+c_n'}{2}$
et ${\d}_{i,n} \leq \frac{2^n+c_n'}{4}$ pour $i \neq \frac{n}{2}$. 
\end{enumerate}
\end{prop}
{\sc Preuve.}
On applique \`a $M=X_m^n$ et $M=Y_{2k}^n$
l'in\'egalit\'e de Smith-Thom
$$b_*({\R}M)=\sum_{i=1}^t b_i({\R}M) \leq b_*(M),$$
o\`u $t =dim \; M$, ainsi que, lorsque $t$ est pair,
les in\'egalit\'es de Comessatti-Petrovsky-Oleinik g\'en\'eralis\'ees
$$2-h^{\frac{t}{2},\frac{t}{2}}(M)
\leq \chi({\R}M)= \sum_{i=1}^{t} {(-1)}^i b_i({\R}M)
\leq h^{\frac{t}{2},\frac{t}{2}}(M).$$
Le r\'esultat d\'ecoule alors de $b_{t-i}({\R}M)=b_i({\R}M) \geq 0$, $b_*(X_m^n)=m^n+{\mathcal R}(n-1,m)$
et de $b_*(Y^n_{2k})=2^n k^n + {\mathcal R}(n-1,k)$.
{\cqfd}
Les bornes ci-dessus sont en g\'en\'eral les meilleures bornes
sup\'erieures connues pour les ${\h}_{i,n}$ et ${\d}_{i,n}$.
Pour \^etre complet, les valeurs de $c_n$ et $c_n'$
sont donn\'ees par les formules suivantes d\'eriv\'ees
de \cite{Hodge}:
$$c_n=\frac{{(-1)}^{\frac{n+1}{2}}}{n!}\sum_{j=1}^{\frac{n+1}{2}}{(-1)}^j
\begin{pmatrix} n+1 \\ \frac{n+1}{2}-j \end{pmatrix}j^n,$$

$$c'_n=\frac{{(-1)}^{\frac{n}{2}+1}}{n!}\left(
\sum_{j=1}^{\frac{n}{2}+1}{(-1)}^j\begin{pmatrix} n+2 \\ \frac{n}{2}+1-j \end{pmatrix}
a_{j,n}+2^n\sum_{j=1}^{\frac{n}{2}}
{(-1)}^j\begin{pmatrix} n+1 \\ \frac{n}{2}-j \end{pmatrix}j^n \right),$$
o\`u $a_{j,n}=\sum_{t=1}^{2j-1}t^n$ et $\begin{pmatrix}a\\b\end{pmatrix}=
\frac{a!}{b!{(a-b)}!}$ est le coefficient binomial.
En particulier, on obtient
$$c_3=\frac{2}{3}, \quad c_5=\frac{11}{20}, \quad c_7=\frac{151}{315},$$
$$c_2'=3, \quad c_4'=\frac{115}{12},\quad c_6'=\frac{5887}{180}.$$

On se propose maintenant d'utiliser notre construction d'hypersurfaces doubl\'ees r\'eelles
afin d'obtenir des bornes inf\'erieures pour les ${\h}_{i,n}$ et ${\d}_{i,n}$.

\begin{prop}
\label{Bet}
Soit $X$ une hypersurface doubl\'ee r\'eelle non singuli\`ere de degr\'e $2k$ dans ${\C}P^n$
dont la partie r\'eelle est hom\'eomorphe au collage
des $${\R}X_k^l(\star) \times {\R}Y_{2k}^{n-l}(\star), \quad l=1,\cdots,n,$$
comme dans la proposition \ref{doub}. Pour tout entier positif $i$, on a
$$b_i({\R}X)=\left(\sum_{l=1}^n \sum_{p=0}^i b_p({\R}X_k^l) \cdot
b_{i-p}({\R}Y_{2k}^{n-l}) \right) +{\mathcal R}(n-1,k).$$
\end{prop}
{\sc Preuve.}
L'in\'egalit\'e de Smith-Thom implique que,
pour tout entier positif $i$, l'on a
$b_i({\R}X_k^l(\star))=b_i({\R}X_k^l)+{\mathcal R}(l-1,k)$ et
$b_i({\R}Y_{2k}^{n-l}(\star))=b_i({\R}Y_{2k}^{n-l})+{\mathcal R}(n-l,k)$.
En utilisant la formule de Kunneth, on obtient ensuite que
$b_i({\R}X_k^l(\star) \times {\R}Y_{2k}^{n-l}(\star))=\sum_{p=0}^i b_p({\R}X_k^l) \cdot
b_{i-p}({\R}Y_{2k}^{n-l})+{\mathcal R}(n-1,k)$.
Le collage de la proposition se faisant le long d'un nombre fini (pour $n$ fix\'e) de vari\'et\'es
de m\^eme type mais de dimension strictement
inf\'erieures \`a $n-1$, l'in\'egalit\'e de Smith-Thom implique que la contribution du collage \` a  $b_i({\R}X)$
est une fonction ${\mathcal R}(n-1,k)$.
{\cqfd}

\begin{prop}
\label{cle}
Supposons donn\'es pour $l=1,\cdots,n$ un polyn\^ome $f_k^l$
non d\'eg\'en\'er\'e de polytope de Newton $T_k^l$ et
un polyn\^ome $f_{2k}^{n-l}$
non d\'eg\'en\'er\'e de polytope de Newton $T_{2k}^{n-l}$.
Il existe une hypersurface r\'eelle $X$ non singuli\`ere de degr\'e $2(k+n+1)$
dans ${\C}P^n$ telle que, pour tout entier positif $i$, l'on ait
$$b_i({\R}X) \geq \left(\sum_{l=1}^n \sum_{p=0}^i b_p({\R}X_k^l) \cdot
b_{i-p}({\R}Y_{2k}^{n-l}) \right) +{\mathcal R}(n-1,k),$$
o\`u $X_k^l$ est l'hypersurface d\'efinie par $f_k^l$ dans ${\C}P^l$ et $Y_{2k}^{n-l}$ est le plan doubl\'e
associ\'e \`a $f_{2k}^{n-l}$ dans ${\C}P^{n-l+1}(1,k)$.
\end{prop}
{\sc Preuve.}
On pose $K=k+n+1$ et on reprend la construction, et donc les notations, de la section \ref{cons} afin d'obtenir
une hypersurface doubl\'ee de degr\'e $2K$ dans ${\C}P^n$
(le $k$ de la section \ref{cons} correspond donc au $K$ ici).

Pour tout entier $l$ compris entre $1$ et $n$, on consid\`ere l'image
du polyn\^ome $g_k^l$ d\'efini par $g_k^l({z}_1,\cdots,{z}_l)= {z}_1 \cdots {z}_lf_k^l({z}_1,\cdots,{z}_l)$
par une transformation affine envoyant $T_K^l$ sur la face $F_1^l$ de $T_K^l$.
On obtient ainsi une collection de polyn\^omes $\tilde{g}_k^l$, $l=1,\cdots,n$, v\'erifiant:
\begin{itemize}
\item $\tilde{g}_k^l$ et $f_k^l$ d\'efinissent la m\^eme hypersurface $X_k^l \cap {({\R}^*)}^l$
dans ${({\R}^*)}^l$,
\item le polytope de Newton $G_k^l$ de $\tilde{g}_k^l$ est contenu \`a l'int\'erieur de $F_1^l$.
\end{itemize}
Il existe alors une subdivision poly\'edrale convexe de $T_K^n$
contenant les $G_k^l$ pour $l=1,\cdots,n$ ainsi qu'un polyn\^ome, de polytope de newton $T_K^n$, dont le tronqu\'e
sur chacun des $G_k^l$ coincide avec le polyn\^ome  $\tilde{g}_k^l$ et dont le tronqu\'e sur tout polytope de la subdivision
est un polyn\^ome non d\'eg\'en\'er\'e. En utilisant la m\'ethode de Viro \`a partir de ces donn\'ees,
on obtient un polyn\^ome $f_K^n$ non d\'eg\'en\'er\'e de polytope de Newton $T_K^n$. L'hypersurface
$X_K^n$ d\'efinie par ce polyn\^ome v\'erifiera alors, pour tout entier positif $i$ et tout entier $l$
compris entre $1$ et $n$
\begin{equation}
\label{1}
b_i({\R}X_K^n \cap X(F_1^l)) \geq b_i({\R}X_k^l)+{\mathcal R}(l-1,k),
\end{equation}
puisque, pour tout entier $l$
compris entre $1$ et $n$, ${\R}X_K^n \cap X(F_1^l)$ contient une copie hom\'eomorphe
de $X_k^l \cap {({\R}^*)}^l$.

En proc\'edant de la m\^eme mani\`ere \`a partir des polyn\^omes $f_{2k}^{n-l}$,
on obtient un polyn\^ome $f_{2K}^n$ de polytope de Newton $T_{2K}^n$ avec la propri\'et\'e suivante.
Si $f_{2K}^{n-l}$ d\'esigne le tronqu\'e de  $f_{2K}^n$ sur $F_2^{n-l}$, alors, pour tout
entier $l$ compris entre $1$ et $n$, l'ensemble $\{{\bar{f}}_{2K}^{n-l} \geq 0\} \subset {\R}X(F_2^{n-l})$ contient
une copie hom\'eomorphe de $\{{\bar{f}}_{2k}^{n-l}\geq 0\} \subset {\R}P^{n-l}$.
Par suite, si $Y_{2K}^{n-l}$ d\'esigne le plan doubl\'e associ\'e \`a $f_{2K}^{n-l}$,
alors, pour tout entier positif $i$ et tout entier $l$ compris entre $1$ et $n$, on a
\begin{equation}
\label{2}
b_i({\R}Y_{2K}^{n-l}) \geq b_i({\R}Y_{2k}^{n-l})+{\mathcal R}(n-l,k).
\end{equation}

Il reste \`a appliquer la proposition \ref{Bet} \`a l'hypersurface doubl\'ee $X$ de polytope de Newton $T_{2K}^n$
construite \` a partir de $f_K^n$ et de $f_{2K}^n$ comme dans la section \ref{cons}.
Le r\'esultat final d\'ecoule ensuite des in\'egalit\'es \ref{1} et \ref{2}.
{\cqfd}

On peut maintenant \'enoncer les r\'esultats principaux de cette section.
\begin{thm}
\label{mth}
$$ {\h}_{0,n} \geq \frac{1}{2^n-2}
\sum_{l=1}^{n-1} {\h}_{0,l} \cdot {\d}_{0,n-l}.$$
\end{thm}
{\sc Preuve.}
La proposition \ref{cle} (pour $i=0$) permet de montrer que
$$Max \; b_0({\R}X^n_{2(k+n+1)}) \geq \sum_{l=1}^n \left( Max \; b_0({\R}X_k^l) \right) \cdot
\left( Max \; b_0({\R}Y_{2k}^{n-l})\right) +{\mathcal R}(n-1,k).$$
En divisant chacun des membres de cette in\'egalit\'e par ${(2k)}^n$ et en faisant tendre
$k$ vers $+\infty$, on obtient
$${\h}_{0,n} \geq \frac{1}{2^n} \sum_{l=1}^{n} {\h}_{0,l} \cdot {\d}_{0,n-l}.$$
Il reste \`a remarquer
que le terme donn\'e par $l=n$ dans la somme vaut ${\d}_{0,0} \cdot {\h}_{0,n}
=2 \cdot {\h}_{0,n}$. {\cqfd}

Soit ${\bar{f}}$ un polyn\^ome r\'eel homog\`ene de degr\'e pair en $n$ variables
d\'efinissant une hypersurface r\'eelle $X$ non singuli\`ere dans ${\C}P^n$.
Une composante connexe de ${\R}X$ est dite positive (resp. n\'egative) si elle borde par l'ext\'erieur
une composante connexe de $\{{\bar{f}} \geq 0\} \subset {\R}P^n$ (resp. $\{{\bar{f}} \leq 0\} \subset {\R}P^n$).
Si l'on suppose que ${\bar{f}}$ est n\'egatif \` a l'ext\'erieur de toute composante connexe de ${\R}X$, alors
le nombre de composantes connexes positives (resp. n\'egatives)
de ${\R}X$ coincide alors avec le nombre de composantes connexes (resp. le nombre de composantes connexes moins $1$)
de la partie r\'eelle du plan doubl\'e
associ\'e \`a $f$ (resp. $-f$).

\begin{thm}
\label{mth1}
$$ {\d}_{0,n} \geq  \sum_{l=1}^{n-1} {\h}_{0,l} \cdot {\d}_{0,n-l} \; + {\h}_{0,n}.$$
\end{thm}
{\sc Preuve.}
On revient \`a la construction
de l'hypersurface doubl\'ee $X$ de la section \ref{cons} (Proposition \ref{doub}).
On remarque que les composantes connexes de ${\R}X$ sont de deux types, celles qui proviennent
de composantes de ${\R}X_{k,+}$ de bords non vides, et celles, venant par paires
(une composante englobe l'autre), provenant des composantes connexes de ${\R}X_k$
contenues dans $\{{\bar{f}}_{2,t} >0\}$. On remarque que les composantes connexes
du premier type sont toutes de m\^eme signe, tandis que dans chaque paire
de composantes connexes du deuxi\`eme type, une composante est positive et l'autre est
n\'egative. On obtient alors de la m\^eme mani\`ere que pour la proposition \ref{Bet}
\begin{equation}
\label{d}
b_0({\R}Y_{2k}^n)= \left(\sum_{l=1}^{n-1} b_0({\R}X_k^l) \cdot
b_0({\R}Y_{2k}^{n-l})\right) + b_0({\R}X_k^n)+{\mathcal R}(n-1,k),
\end{equation}
si $Y_{2k}^n$ est l'un des plans doubl\'es au dessus de $X$
i.e. $Y=\{U^2 \pm {\bar{f}}_{2k}=0\} \subset {\C}P^{n+1}(1,k)$ si $X=\{{\bar{f}}_{2k}=0\}$.
Ensuite, comme dans la preuve de la proposition \ref{cle}, on montre que dans l'\'egalit\'e \ref{d}, quitte \`a remplacer
${\R}Y_{2k}^n$ par ${\R}Y_{2(k+c)}^n$ avec $c$ ne d\'ependant pas de $k$, on peut supposer
que les $X_k^l$ et $Y_{2k}^{n-l}$ sont arbitraires. Finalement, en prenant les maxima de chacun
des nombres de betti de cette in\'egalit\'e, en divisant par $k^n$, puis en faisant tendre $k$ vers $+\infty$
(voir la preuve du th\'eor\`eme \ref{mth}), on obtient le r\'esultat souhait\'e.
{\cqfd}
\smallskip

On se concentre maintenant sur les surfaces alg\'ebriques r\'eelles
dans ${\C}P^3$.

\begin{thm}
\label{surfg}

$$ \frac{{\d}_{0,2}}{6}+\frac{1}{12} \leq {\h}_{0,3} \leq \frac{5}{12}, \quad \frac{{\d}_{1,2}}{6}+\frac{1}{6}
\leq {\h}_{1,3} \leq \frac{5}{6}.$$
\end{thm}
{\sc Preuve.}
Les bornes sup\'erieures sont celles de la proposition
\ref{bsup}. L'in\'egalit\'e
${\h}_{0,3} \geq \frac{{\d}_{0,2}}{6}+\frac{1}{12}$ s'obtient \`a partir du th\'eor\`eme \ref{mth} en utilisant les valeurs connues
${\h}_{0,1}={\d}_{0,1}=1$ et ${\h}_{0,2}=1/2$. Montrons maintenant l'in\'egalit\'e ${\h}_{1,3} \geq \frac{{\d}_{1,2}}{6}+\frac{1}{6}$.
Pour $n=3$ et $i=0$, l'in\'egalit\'e de la proposition \ref{cle} donne (apr\`es simplifications,
certains nombres de Betti sont nuls pour des raisons de dimension, d'autres sont \'egaux deux \`a deux par dualit\'e):
$$
\begin{array}{lll}
b_1({\R}X) & \geq & b_0({\R}X_k^1) \, \cdot \,
b_{1}({\R}Y_{2k}^{2}) \; + \; 2b_0({\R}X_k^2) \, \cdot \,
b_{1}({\R}Y_{2k}^{1})\; + \\
 & & \\
&  & b_1({\R}X_k^3) \, \cdot \,
b_{0}({\R}Y_{2k}^{0}) \; + \;{\mathcal R}(2,k)
\end{array}
.$$
En prenant les maxima de chacun des nombres de Betti
de cette derni\`ere in\'egalit\'e, en divisiant les deux membres par $(2k)^3$,
puis en faisant tendre $k$ vers $+\infty$, on obtient
${\h}_{1,3} \geq \frac{1}{8}({\h}_{0,1} \cdot {\d}_{1,2}+2{\h}_{0,2} \cdot {\d}_{1,1}+{\h}_{1,3} \cdot {\d}_{0,0})$.
En utilisant les valeurs connues des termes du membre de droite, on obtient le r\'esultat.
{\cqfd}

\begin{cor}
\label{surf}
$$ \frac{35}{96} \leq {\h}_{0,3} \leq \frac{5}{12}, \quad \frac{35}{48}
\leq {\h}_{1,3} \leq \frac{5}{6}.$$
\end{cor}
{\sc Preuve.}
C'est une cons\'equence du th\'eor\`eme pr\'ec\'edent ainsi que des in\'egalit\'es
${\d}_{0,2} \geq 27/16 \; , \; {\d}_{1,2} \geq \frac{27}{8}$ qui sont prouv\'ees
dans \cite{Itrag}.
{\cqfd}

\begin{rem}
\label{gen}
Les in\'egalit\'es $\frac{13}{36} \leq {\h}_{0,3}$ et $\frac{13}{18} \leq {\h}_{1,3}$
de \cite{Basymp} sont obtenues de mani\`ere similaire \`a partir
des in\'egalit\'es ${\d}_{0,2} \geq \frac{5}{3}$ et
${\d}_{1,2} \geq \frac{10}{3}$ obtenues dans \cite{Ha}. En fait, on peut extraire de \cite{Basymp}
(voir \cite{Basymp}, remarque 2, section 5) les in\'egalit\'es
${\h}_{0,3} \geq \frac{T{\d}_{0,2}}{6}+\frac{1}{12}$ et ${\h}_{1,3} \geq \frac{T{\d}_{1,2}}{6}+\frac{1}{6}$
o\`u $T{\d}_{0,2}$ et $T{\d}_{1,2}$ sont d\'efinis comme ${\d}_{0,2}$ et ${\d}_{1,2}$ mais en se restreignant
aux $T$-courbes (voir la sous-section \ref{T-hyp} pour une d\'efinition pr\'ecise).
\end{rem}

On pr\'esente maintenant un tableau recensant les bornes inf\'erieures pour les ${\h}_{0,n}$
et ${\d}_{0,n}$, $n \leq 7$, obtenues r\'ecursivement grace au th\'eor\`eme \ref{mth1}
\`a partir des estimations connues pour $n \leq 2$.
Les bornes sup\'erieures sont celles de la proposition \ref{bsup}
sauf pour les ${\d}_{0,n}$ avec $n$ impair pour lesquels
on a utilis\'e l'in\'egalit\'e ${\d}_{0,n} \leq 2^n {\h}_{0,n}$ (voir la remarque \ref{rqb}
plus bas) alli\'ee avec la borne sup\'erieure pour ${\h}_{0,n}$ de la proposition \ref{bsup}.
\bigskip

\begin{center}
\begin{tabular}{|c|c|c|c|c|} \hline

$n$ & \multicolumn{2}{c|} {${\h}_{0,n}$} & \multicolumn{2}{c|} {${\d}_{0,n}$} \\[3mm] \hline
 & inf & sup & inf & sup \\[3mm] \hline
$1$ & $\;1\;$&$\;1\;$ &$\;1\;$ &$\;1$\; \\[3mm] \hline
& & & & \\
$2$ & $\frac{1}{2}$& $\frac{1}{2}$ & $\frac{27}{16}$& $\frac{7}{4}$ \\[3mm] \hline
& & & & \\
$3$ & $\frac{35}{96}$& $\frac{5}{12}$& $\frac{245}{96}$& $\frac{10}{3}$ \\[3mm] \hline
& & & & \\
$4$ & $\frac{361}{1344} \sim 0,27$& $\frac{1}{2}$& $\frac{1805}{448} \sim 4,03$& $\frac{307}{48} \sim
6,4$\\[3mm] \hline
& & & & \\
$5$ & $\frac{22181}{107520} \sim 0,2$& $\frac{31}{80}$& $\frac{687611}{107520} \sim 6,4$& $\frac{62}{5}$\\[3mm]
\hline
& & & & \\
$6$ & $\;\frac{1612753}{9999360} \sim 0,16\;$& $\;\frac{1}{2}\;$& $\;\frac{1612753}{158720} \sim 10,16\;$&
\;$\frac{17407}{720} \sim 24,17\;$
\\[3mm] \hline
& & & & \\
$7$ & $\;\frac{854473649}{6719569920} \sim 0,127\;$& $\;\frac{233}{630} \sim 0,37\;$& $\;
\frac{108518153423}{6719569920} \sim 16,15\;$
& $\; \frac{14912}{315} \sim 47,34\;$\\[3mm] \hline
\end{tabular}
\end{center}
\medskip

\begin{center}
{\sc Tableau 1.}
\end{center}
\bigskip

\begin{rem}
\label{rqb}
On a $$2^{n-1}{\h}_{0,n} \leq {\d}_{0,n} \leq 2^n {\h}_{0,n}$$
car une composante connexe d'une hypersurface projective
${\R}X_{2k}^n=\{{\bar{f}}_{2k}=0\} \subset {\R}P^n$
est soit positive soit n\'egative, ce qui implique que
$\frac{1}{2}b_0({\R}X_{2k}^n) \leq b_0({\R}Y_{2k}^n) \leq b_0({\R}X_{2k}^n)$
pour l'un ou l'autre des plans doubl\'es $Y_{2k}^n=\{U^2 \pm {\bar{f}}_{2k}(Z)=0\}$.
\end{rem}
On peut trouver une borne moins bonne que celle du th\'eor\`eme \ref{mth1}, mais beaucoup plus facilement
calculable.

\begin{prop}
\label{le}
$${\h}_{0,n} \geq \frac{1}{2^{n-1}}.$$
\end{prop}
{\sc Preuve.}
On obtient le r\'esultat par r\'ecurrence \`a partir des in\'egalit\'es
${\h}_{0,n} \geq \frac{1}{2^n-2} \sum_{l=1}^{n-1} {\h}_{0,l} \cdot {\d}_{0,n-l}$ et
${\d}_{0,n} \geq 2^{n-1} {\h}_{0,n}$ du th\'eor\`eme \ref{mth1} et de la remarque \ref{rqb}.
{\cqfd}
\smallskip

\subsection{Cas des T-hypersurfaces}
\label{T-hyp}

Un $T$-polyn\^ome est un polyn\^ome de Viro $f_t$ pour lequel la subdivision poly\'edrale (de son polytope de Newton)
associ\'ee est une triangulation
dont l'ensemble des sommets correspond exactement \`a l'ensemble des mon\^omes de $f_t$. Une $T$-hypersurface
est une hypersurface d\'efinie, pour $t>0$ suffisamment petit, par un $T$-polyn\^ome. Notons que g\'en\'eriquement
un polyn\^ome de Viro est un $T$-polyn\^ome. Le type topologique de la partie r\'eelle d'une $T$-hypersurface
d\'efinie par un $T$-polyn\^ome $f_t$ ne d\'epend que des signes (et donc pas de la valeur absolue) des coefficients de $f_t$
et de la triangulation associ\'ee. La proc\'edure permettant d'obtenir le type topologique de la partie r\'eelle d'une $T$-hypersurface
\`a partir de la triangulation associ\'ee et de la distribution de signes aux sommets de cette triangulation
est appel\'ee patchwork combinatoire ou $T$-construction (pour une description pr\'ecise du patchwork combinatoire,
voir, par exemple, \cite{ItSh,ItVi}). La famille des $T$-hypersurfaces de degr\'e et dimension donn\'es est assez rigide,
(on peut noter, par exemple, que si $f$ est un $T$-polyn\^ome affine en $n$ variables et $a$ est un changement de coordonn\'ees affines
, alors en g\'en\'eral, le polyn\^ome $f \circ a$ n'est pas un $T$-polyn\^ome) et l'on doit
s'attendre \`a ce qu'elle soit strictement contenue
dans la famille correspondante d'hypersurfaces alg\'ebriques r\'eelles. Itenberg \cite{ItT-cour} a montr\'e l'existence de
courbes alg\'ebriques r\'eelles planes qui ne sont pas des $T$-courbes. Plus r\'ecemment, Itenberg et Shustin \cite{ItSh} ont montr\'e
l'existence d'hypersurfaces alg\'ebriques r\'eelles dans ${\C}P^n$ qui ne sont pas des $T$-hypersurfaces
pour tout $n \geq 7$ et tout $m$ suffisamment grand. Le but ici est d'obtenir le r\'esultat similaire
pour des dimensions inf\'erieures.

On commence par remarquer que la proposition \ref{as1} a son \'equivalent direct pour les $T$-hypersurfaces.

Soit $i$ un entier positif. On note, respectivement,
$$TMax \; b_i({\R}  X_m^n) \quad \mbox{et} \quad TMax \; b_i({\R}Y_{2k}^n) $$
la valeur maximale
des nombres de Betti $b_i({\R}  X_m^n)$ prise sur l'ensemble des $T$-hypersurfaces $X_m^n$
de degr\'e $m$ dans ${\C}P^n$ pour $m$ et $n$ fix\'es, et la valeur maximale des nombres de Betti $b_i({\R}Y_{2k}^n)$
prise sur l'ensemble des plans doubl\'es $Y_{2k}^n=\{U^2-{\bar{f}}_{2k}(Z)=0\} \subset {\C}P^{n+1}(1,k)$ associ\'es
\`a des des $T$-polyn\^omes $f_{2k}$ pour $k$ et $n$ fix\'es.

\begin{prop}
\label{as1T}
Il existe des nombres r\'eels $T{\h}_{i,n}$ et $T{\d}_{i,n}$
pour lesquels on ait les \'equivalences asymptotiques suivantes
$$ TMax \; b_i({\R}  X_m^{n}) \sim T{\h}_{i,n} \cdot m^n, \quad
TMax \; b_i({\R}Y_{2k}^{n}) \sim T{\d}_{i,n} \cdot k^n$$
pour $m$ (resp. $k$) $\rightarrow + \infty$.
\end{prop}

{\sc preuve.}
On montre exactement comme dans la preuve de la proposition \ref{as1}
que la suite $TMax \; b_i({\R}X_m^n)/m^n$ converge vers sa limite sup\'erieure $L$.
On se contente d'indiquer ici les modifications \` a apporter de mani\`ere \`a ce que
les hypersurfaces construites dans la preuve de la proposition \ref{as1} soient des $T$-hypersurfaces
v\'erifiant les m\^emes propri\'et\'es.

Soit  $X_{m_0}^n$ une $T$-hypersurface telle que
$b_i({\R}X_{m_0}^n)/{m_0}^n>L-\epsilon/2$. On note $\tau$ la triangulation convexe
de $T_{m_0}^n$ et $D_{\tau}$ la distribution de signes aux sommets de $\tau$ qui sont associ\'ees \`a $X_{m_0}^n$.

Par la suite, l'image de $\tau$ par une transformation affine $\Delta\in AFF_n(\Z)$
sera la triangulation convexe $\Delta(\tau)$ de $\Delta(T_{m_0}^n)$
constitu\'ee des $\Delta(\gamma)$ pour $\gamma \in \tau$ et l'image de $D_{\tau}$ par $\Delta$ sera la distribution de
signes $\Delta(D_{\tau})$ aux sommets de $\Delta(\tau)$ attribuant \`a $\Delta(w)$ le signe attribu\'e par $D_{\tau}$ \`a $w$.

On peut supposer, quitte \`a augmenter $m_0$,
que la $T$-hypersurface $X_{m_0}^n$ n'intersecte pas les hyperplans de coordonn\'ees de ${\R}P^n$.
En effet, soit $\Delta \in AFF_n(\Z)$ la translation par $(1,\cdots,1)$ et soit $c=1$ ou $2$ tel que
${m_0+n+c}$ est pair. On \'etend $\Delta(\tau)$ en une triangulation convexe $\tilde{\tau}$ de $T_{m_0+n+c}^n$
de telle sorte que chacun des sommets de $\tilde{\tau}$ situ\'e sur le bord de $T_{m_0+n+c}^n$ appartienne \`a ${(2 \Z)}^n$,
puis on \'etend $\Delta(D_{\tau})$ en une distribution de signes $D_{\tilde{\tau}}$ de telle sorte que $D_{\tilde{\tau}}$ attribue
le m\^eme signe \`a tous les sommets situ\'e sur le bord de $T_{m_0+n+c}^n$.
Clairement, si $X_{m_0+n+c}^n$ est une $T$-hypersurface associ\'ee \`a $\tilde{\tau}$ et
$D(\tilde{\tau})$, alors $X_{m_0+n+c}^n$ n'intersecte pas les hyperplans de coordonn\'ees de ${\R}P^n$ et v\'erifie
$b_i({\R}X_{m_0+n+c}^n)/{m_0}^n>L-\epsilon/2$ pour $m_0$ suffisamment grand.

Pour construire les hypersurfaces $X_{p({m_0}+n+1)}^n$ comme des $T$-hypersurfaces, on
utilise une triangulation convexe de $T_{p({m_0}+n+1)}^n$ ainsi qu'une distribution de signes
qui \'etendent les $p^n$ images correspondantes de $\tau$ et de $D(\tau)$.

Finalement, pour compl\'eter la famille de $T$-hypersurfaces 
$X_{p({m_0}+n+1)}^n$, $p \geq p_0$, en une famille de $T$-hypersurfaces ${\mathcal F}=\{X_m^n, \; m \geq p_0({m_0}+n+1) \}$,
on consid\`ere, pour tout entier $m$ tel que $p({m_0}+n+1)< m < (p+1)({m_0}+n+1)$ avec $p \geq p_0$, 
une triangulation convexe de $T_m^n$ et une distribution de signes \`a ses sommets
qui \'etendent la triangulation convexe de $T_{p({m_0}+n+1)}^n$ et distribution de signes correspondante.

La preuve de l'existence de $T\d_{i,n}$ s'obtient \`a partir de la preuve de l'existence de $\d_{i,n}$
par les modifications similaires.
{\cqfd}

On a bien sur $T\h_{i,n} \leq h_{i,n}$ et $T\d_{i,n} \leq d_{i,n}$. Pour les petites dimensions,
on a en fait des \'egalit\'es:

$$
{\d}_{0,0}=T{\d}_{0,0}=2,\quad
{\h}_{0,1}={\d}_{0,1}={\d}_{1,1}=T{\h}_{0,1}=T{\d}_{0,1}=T{\d}_{1,1}=1,$$
$$
{\h}_{0,2}={\h}_{1,2}=T{\h}_{0,2}=T{\h}_{1,2}=\frac{1}{2} \quad (\mbox{voir, par exemple, \cite{ItT-cour}}).
$$

On connait de plus les estimations suivantes \cite{Itrag}
$$T{\d}_{0,2} \geq \frac{27}{16}, \quad T{\d}_{1,2} \geq \frac{27}{8}$$
ainsi que
$$
\frac{1}{4} \leq T{\h}_{0,3} \leq \frac{16}{39}, \quad \frac{97}{144} \leq T{\h}_{1,3} \leq \frac{7}{9}.$$
Les bornes inf\'erieures pour $T{\h}_{0,3}$ et $T{\h}_{1,3}$ proviennent de \cite{these} et \cite{T-surf}, respectivement,
les bornes sup\'erieures proviennent d'un article r\'ecent \cite{ItSh} d'Itenberg et Shustin, duquel est extrait
la proposition suivante.

\begin{prop}[\cite{ItSh}]
\label{propI}
Pour tout $n \geq 4$, on a
$$T{\h}_{0,n} \leq \frac{2^{n-1}}{n!}.$$
\end{prop}

Comme cons\'equence de cette proposition et de notre construction, on obtient le r\'esultat suivant.

\begin{thm}
\label{Tsurf}
Pour tout $n \geq 5$, on a
$$T{\h}_{0,n} < {\h}_{0,n}.$$
En particulier, pour tout $n \geq 5$ et tout degr\'e $m$ suffisamment grand,
il existe des hypersurfaces alg\'ebriques r\'eelles $X_m^n$ de degr\'e $m$ dans ${\C}P^n$
qui ne sont pas des $T$-hypersurfaces.
\end{thm}

{\sc Preuve.}
On obtient $n \geq 5 \Rightarrow {\h}_{0,n} > \frac{2^{n-1}}{n!}$
grace au lemme \ref{le} pour $n \geq 7$ et du tableau 1.
pour $n=5$ et $6$. Il reste \`a appliquer la proposition \ref{propI}.
{\cqfd}
\bigskip

{\bf Remarques conclusives.}
\medskip

{\bf 1.} Une vari\'et\'e alg\'ebrique r\'eelle lisse $X$
est appel\'ee {\it M-vari\'et\'e} si elle rend exacte l'in\'egalit\'e
de Smith-Thom i.e. si $b_*({\R}X)=b_*(X)$.
Il existe des M-hypersurfaces $X_m^n$ pour tout $m$ et $n$ \cite{ItVi}.
Une famille d'hypersurfaces $X_m^n$ est dite {\it asymptotiquement maximale}
si $b_*({\R}X_m^n)=b_*(X_m^n)+{\mathcal R}(n-1,m)=m^n+{\mathcal R}(n-1,m)$.
On remarque que la famille des hypersurfaces doubl\'ees construites
dans la section \ref{cons} est asymptotiquement maximale
si c'est la cas pour chacune des familles d'hypersurfaces $X_k^l$ et $X_{2k}^{n-l}$.
\medskip

{\bf 2.} On peut montrer des r\'esultats du m\^eme type que la proposition \ref{as1}
pour d'autres invariants topologiques que les nombres de Betti individuels.
Par exemple, en reprenant la preuve de la proposition \ref{as1}, on montre ais\'ement
l'existence d'un r\'eel $\chi_n$ tel que $Max \; \chi({\R}X_m^n) \sim \chi_n \cdot m^n$.

\bibliographystyle{amsalpha}

\end{document}